\documentclass[journal]{IEEEtran}

\usepackage{amsmath,amssymb}
\usepackage{xcolor}
\usepackage{amsfonts}
\usepackage{bm}
\usepackage{graphicx}
\usepackage{multicol}
\usepackage{subfigure}
\usepackage{float}
\usepackage{epsfig}
\usepackage{epstopdf}
\usepackage{pdfpages}
\usepackage{booktabs}
\usepackage{mathtools}
\usepackage[linesnumbered,ruled,vlined]{algorithm2e}

\newtheorem{lemma}{Lemma}

\newtheorem{theorem}{Theorem}

\newcommand{\comment}[1]{}

\renewcommand\IEEEkeywordsname{Keywords}

\DeclarePairedDelimiter\abs{\lvert}{\rvert}%

\newcommand{\nosemic}{\renewcommand{\@endalgocfline}{\relax}}% Drop semi-colon ;
\newcommand{\dosemic}{\renewcommand{\@endalgocfline}{\algocf@endline}}% Reinstate semi-colon ;
\let\oldnl\nl% Store \nl in \oldnl
\newcommand{\nonl}{\renewcommand{\nl}{\let\nl\oldnl}}% Remove line number for one line

\title{Exact aggregate models for optimal management of heterogeneous fleets of storage devices}
\author{David~Angeli,~\IEEEmembership{Fellow, IEEE},
        Zihang~Dong,~and~Goran~Strbac,~\IEEEmembership{Member, IEEE}% <-this % stops a space

\thanks{
The authors are with the Department of Electrical and Electronic Engineering, Imperial College London, London, SW7 2AZ, UK. {\textit{ (e-mail:  d.angeli@imperial.ac.uk;
zihang.dong14@imperial.ac.uk; g.strbac@imperial.ac.uk)}}}%
\thanks{David Angeli is also with the Dipartimento di Ingegneria dell'Informazione, Universit\`a di Firenze, Italy}
}
\begin{document}
\maketitle
\begin{abstract}
Future power grids will entail large fleets of storage devices capable of scheduling their charging/discharging profiles so as to achieve lower peak demand and reducing energy bills, by shifting absorption times in sync with availability of renewable energy sources. Optimal management of such fleets entails large scale optimisation problems which are better dealt with in a hierarchical manner, by clustering together individual devices into fleets. Leveraging on recent results characterizing the set of aggregate demand profiles of a heterogeneous fleet of charging (or, respectively, discharging) devices we propose a way to achieve optimality, in a unit commitment problem, by adopting a simplified formulation with a number of constraints for the fleet that scales linearly in the number of time-slots considered and is independent of the size of the fleet. This is remarkable, as it shows that, under suitable conditions, a heterogeneous fleet of any size can effectively be treated as a single storage unit.
\end{abstract}


\section{Introduction and Motivations}
Recent years have witnessed a considerable interest in the scientific community towards \emph{optimal management} of distributed energy resources. It is anticipated that implementation of the low-carbon agenda and the adoption of 
intermittent renewable energy sources will exacerbate issues of demand-supply balance and significantly incentivize the adoption of flexible demand paradigms across the grid. 

A recent survey, \cite{weitzel}, has systematically classified a
 significant amount of literature dealing with optimal energy management for Electric Vehicles alone (EVs). 
While a fraction of such papers is concerned only with scheduling aspects related to a single EV and the benefits that this could enable for the household, many authors are also concerned with the issue of how solutions implemented at the level of the single household may impact the overall system or how to suitably manage (and suitably distribute) the resulting
 large scale optimisation problem.

An emerging viewpoint, in such respect, is that a hierarchical approach, \cite{bernstein}, may allow to cluster multiple households, EVs or more generally prosumers, into groups which are regarded as a single entity by higher level grid abstraction layers, so that only their aggregate power absorption or generation profiles need to be considered, without mention of lower-level detailed dynamics and schedules.
The wisdom is that this could provide an effective modular approach for tackling the inherent complexity of the optimisation and scheduling task. 

To make this proposition effective, however, it is important to characterize the degrees of freedom afforded by an aggregated
fleet of flexible devices in a way that is accurate and still
simple. For instance, \cite{hill}, introduces a heuristic aggregated model of flexible demand, to be used for the solution of Unit Commitment problems.
Similarly, \cite{calvillo}, models the aggregate behaviour of \emph{inelastic} flexible demand through heuristic constraints that limit the shifting ability afforded throughout the day. 
Promising results, in this respect, are provided in \cite{ETA-pscc, ETA-tops, ETA-tosg, EAST} where optimal causal scheduling policies (affording maximal flexibility) for heterogeneous fleets and arbitrary time-varying demand signals are provided and their intrinsic capabilities characterized graphically through so called E-P transforms. Equivalent time-domain characterizations were independently explored in \cite{zachary} and further expanded and interpreted in \cite{phtr}. 
A related line of investigation, \cite{appino}, proposes the use of controlled invariant sets in order to characterize viable power profiles for aggregated fleets, so as to reduce their model complexity for the purpose of higher level optimisation.

More recently, the challenging task of designing optimal dispatch policies for devices with \emph{partial availability} (as typical of Electric Vehicles which are not always connected to the grid) was considered in \cite{submitted,cdcbatteries}.
In particular, \cite{submitted} proposes a scheduling algorithm for fleets with partial availability and at the same time characterizes the set of aggregated demand profiles that  heterogeneous fleets with partial flexibility induce, when operated without \emph{cross-charging}. 
For a fleet with $N$ storage devices and a horizon of $T$ time
intervals, 
this replaces a description involving $N \cdot T$ variables and $2 \cdot N \cdot T$ constraints with one that only entails $T$ variables (the aggregated demand signal) and
$2^T$ constraints, thus drastically improving the scalability with respect to $N$, the size of the fleet.
In this paper we propose a way to achieve optimality, in a unit commitment problem, by adopting a simplified formulation with a number of constraints for the fleet that scales linearly in the number $T$ of time-slots considered and is also independent of the size of the fleet $N$. This is remarkable, as it shows that, under suitable conditions, a heterogeneous fleet of any size can effectively be treated as a single storage unit. In addition, we provide the relevant theory in support of our simplification algorithms and compare the different approaches on a large number of randomly generated examples and case studies.

The remainder of this paper is organized as follows: Section \ref{sec:aggregatemodel} presents the aggregate model of the storage fleet. Section \ref{sec:fullavailability} introduces the optimization problem for devices with full availability window. The aggregated model and the associated optimization algorithms for partial available storage devices are elaborated in \ref{sec:partialavailability}. The advantages of proposed aggregated model in computation is addressed in \ref{sec:computation}. Section \ref{sec:twoarea} applies the algorithms to a two-area system, and Section \ref{sec:conclusion} concludes this paper.

\section{An aggregate fleet model for optimisation}\label{sec:aggregatemodel}
We consider a finite discrete time horizon $\mathcal{T} = \{1,2,\ldots,T\}$ and a storage fleet with $N$ batteries, which we denote
$\mathcal{N} = \{1,2, \ldots, N\}$. Each battery $j \in \mathcal{N}$ is initially empty (without loss of generality) and needs to reach a certain target energy level
$\bar{E}_j$ within its availability set $\mathcal{A}_j$.
Moreover, each battery has a maximum rated power $\bar{P}_j$ which is the maximum power it can absorb while charging.
We consider the problem of minimizing generation costs for the task of charging the fleet while additionally meeting some assigned inflexible demand profile 
$D^{I}(t)$. We assume that a set of generation units $\mathcal{I}$ is available to this end, each one with minimum and maximum power ratings
$\underline{G}_i$ and $\bar{G}_i$, respectively, $i \in \mathcal{I}$.
In particular the following optimisation problem is of interest:
\begin{subequations}\label{original}
\begin{align}
\min_{\begin{aligned}
	u_j(t)\,\,& j \in \mathcal{N}, t \in \mathcal{T} \\ g_i(t)\,\,& i \in \mathcal{I}, t \in \mathcal{T} 
\end{aligned} } &  \sum_{t \in \mathcal{T}} \sum_{i \in \mathcal{I} }   C_i (g_i(t)) \tag{\ref{original}}\\
s.t.\quad \sum_{i \in \mathcal{I}} g_i(t) &= D^{I}(t) + \sum_{j \in \mathcal{N}} u_j(t) \quad \forall \, t \in \mathcal{T} \\
0 \leq &u_j(t) \leq \bar{P}_j \quad \forall \, t \in \mathcal{T}, \; \forall \, j \in \mathcal{N}\\
\underline{G}_i \leq & g_i(t) \leq \bar{G}_i \quad \forall \, t \in \mathcal{T}, \; \forall \, i \in \mathcal{I} \\
\bar{E}_j &\leq \sum_{t \in \mathcal{T}} u_j (t) \quad \forall \, j \in \mathcal{N} \\
u_j(t) &= 0 \quad \forall \, t \notin \mathcal{A}_j
\end{align}
\end{subequations}
% \begin{equation}
% \label{original}
% \min_{\begin{array}{rl}
% 	u_j (t)& j \in \mathcal{N}, t \in \mathcal{T} \\ g_i(t) & i \in \mathcal{I}, t \in \mathcal{T} 
% \end{array} }  \sum_{t \in \mathcal{T}} \sum_{i \in \mathcal{I} }   C_i (g_i(t)) 
% \end{equation}
% \[ \textrm{subject to} \]
% \[
% \begin{array}{rl}
% \sum_{i \in \mathcal{I}} g_i(t) &= D^{I}(t) + \sum_{j \in \mathcal{N}} u_j(t) \qquad \forall \, t \in \mathcal{T}  \\
% 0 \leq u_j(t) \leq& \bar{P}_j \qquad \forall \, t \in \mathcal{T}, \; \forall \, j \in \mathcal{N} \\ 
% \underline{G}_i \leq g_i(t) \leq & \bar{G}_i \qquad \forall \, t \in \mathcal{T}, \; \forall \, i \in \mathcal{I} \\
% \sum_{t \in \mathcal{T}} u_j (t) & \geq \bar{E}_j \qquad \forall \, j \in \mathcal{N} \\
% u_j(t) &= 0 \qquad \forall \, t \notin \mathcal{A}_j
% \end{array}   
% \]	
where $C_i(\cdot)$ are convex increasing functions on their domain
$[\underline{G}_i, \bar{G}_i]$.
Let us consider the set of possible aggregated demand signals of a fleet of batteries operating without cross-charging and within given availability windows:
\begin{equation*}
\begin{aligned}
\mathcal{D}:= \Big \{ d: \mathcal{T} \rightarrow \mathbb{R}:
\exists \, u_j: \mathcal{T} \rightarrow [0, \bar{P}_j]: d(t) = \sum_{j \in \mathcal{N}} u_j(t), & \\
u_j(t) = 0 \; \forall \, j \in \mathcal{N}, \forall \, t \notin \mathcal{A}_j, \text{ and } \sum_{t \in \mathcal{T}} u_j(t) \geq \bar{E}_j,  &
 \Big \}.
\end{aligned}
\end{equation*}
% \[ \mathcal{D}:= \Big \{ d: \mathcal{T} \rightarrow \mathbb{R}:
% \exists \, u_j: \mathcal{T} \rightarrow [0, \bar{P}_j]:
% u_j(t) = 0 \; \forall \, j \in \mathcal{N}, \forall \, t \notin \mathcal{A}_j,
% \]
% \[ \qquad \qquad 
% d(t) = \sum_{j \in \mathcal{N}} u_j(t), \textrm{ and } \sum_{t \in \mathcal{T}} u_j(t) \leq \bar{E}_j,  
%  \Big \}. \]
This polyhedron is described implicitly by using $2NT + N$ inequality constraints and $T+TN$ variables, with (typically) $N \gg T$, (i.e. $T=24$, while $N$ could be in the thousands or even larger). \\
Alternatively, one may seek to describe the set of feasible aggregate power demands without explicit reference to individual battery profiles, viz. by only considering constraints on the aggregate demand. This is possible since the set $\mathcal{D}$ is convex and essentially a projection of a polytope
(of dimension $NT + N$) onto a $T$-dimensional polytope. 
It was shown in \cite{submitted} that $\mathcal{D}$ can also be characterized as follows:
\begin{equation*}
\begin{aligned}
&\mathcal{D} = \Big \{ d: \mathcal{T} \rightarrow \mathbb{R}_{\geq 0}: \\
&\sum_{t \in W} d(t) \leq \sum_{j \in \mathcal{N}} \min \{ \textrm{card} 
( \mathcal{A}_j \cap W) , \bar{E}_j / \bar{P}_j \} \bar{P}_j, \; \forall \, W \subseteq \mathcal{T} \Big \}.
\end{aligned}
\end{equation*}
% \[ \mathcal{D} = \Big \{ d: \mathcal{T} \rightarrow \mathbb{R}_{\geq 0}: 
% \sum_{t \in W} d(t) \leq \sum_{j \in \mathcal{N}} \min \{ \textrm{card} 
% ( \mathcal{A}_j \cap W) , \bar{E}_j / \bar{P}_j \} \bar{P}_j, \; \forall \, W \subseteq \mathcal{T} \Big \}
% \]
As a consequence, the previous optimisation problem can be equivalently recast as:
\begin{subequations}\label{modified}
\begin{align}
\min_{\begin{aligned}
	d (t)\,\,& t \in \mathcal{T} \\ g_i(t)\,\,& i \in \mathcal{I}, t \in \mathcal{T} 
	\end{aligned} }&  \sum_{t \in \mathcal{T}} \sum_{i \in \mathcal{I} }   C_i (g_i(t)) \tag{\ref{modified}}\\
s.t.\quad \sum_{i \in \mathcal{I}} g_i(t) &= D^{I}(t) + d(t) \quad \forall \, t \in \mathcal{T}  \\
\underline{G}_i \leq & g_i(t) \leq \bar{G}_i \quad \forall \, t \in \mathcal{T}, \; \forall \, i \in \mathcal{I} \\
\sum_{t \in \mathcal{T}} d(t) & = \sum_{j \in \mathcal{N}} \bar{E}_j \\
\sum_{t \in \mathcal{W}} d(t) & \leq \sum_{j \in \mathcal{N}}
\min \{ \textrm{card} ( \mathcal{A}_j \cap \mathcal{W}), \bar{E}_j / \bar{P}_j \} \bar{P}_j \nonumber \\
& \qquad \qquad \qquad \qquad \qquad \quad  \forall \, \mathcal{W} \subset \mathcal{T}. 
\end{align}
\end{subequations}

% \begin{equation}
% \label{modified}
% \min_{\begin{array}{rl}
% 	d (t)&  t \in \mathcal{T} \\ g_i(t) & i \in \mathcal{I}, t \in \mathcal{T} 
% 	\end{array} }  \sum_{t \in \mathcal{T}} \sum_{i \in \mathcal{I} }   C_i (g_i(t)) 
% \end{equation} 
% \[ \textrm{subject to} \]
% \[
% \begin{array}{rl}
% \sum_{i \in \mathcal{I}} g_i(t) &= D^{I}(t) + d(t) \qquad \forall \, t \in \mathcal{T}  \\ 
% \underline{G}_i \leq g_i(t) \leq & \bar{G}_i \qquad \forall \, t \in \mathcal{T}, \; \forall \, i \in \mathcal{I} \\
% \sum_{t \in \mathcal{W}} d(t) & \leq \sum_{j \in \mathcal{N}}
% \min \{ \textrm{card} ( \mathcal{A}_j \cap \mathcal{W}), \bar{E}_j / \bar{P}_j \} \bar{P}_j \qquad \forall \, \mathcal{W} \subset \mathcal{T} \\
% \sum_{t \in \mathcal{T}} d(t) & = \sum_{j \in \mathcal{N}} \bar{E}_j 
% \end{array}   
% \]	
\section{A simplified model for full availability}\label{sec:fullavailability}

While the previous model does not involve detailed descriptions of individual batteries (in that it only involves aggregated demand profiles) it still requires a potentially large number of inequality constraints, $2^T-1$, in order to capture the flexibility afforded by the fleet. Its use is therefore beneficial only when $2^T < 2NT+T$, viz. for very large fleets or for a relatively low number of considered time intervals.

We perform in Section \ref{sec:computation} a systematic analysis of how
problem (\ref{modified}) and (\ref{original}) compare from the 
numerical point of view.
While some of the constraints in
(\ref{modified}) will typically be redundant, their a priori elimination appears to be too computationally demanding as it entails the solution of an LP of comparable computational complexity as the problem itself.
We propose next a way to replace the $2^T$ constraints with $T$ inequality constraints which still afford to compute the true optimal solution (or a very close approximation of it). \\

We consider first the case of full availability, viz. $\mathcal{A}_j =
\mathcal{T}$ for all $j \in \mathcal{N}$. Under such scenario,
let $t_1,t_2,\ldots, t_T$ be such that:
\[  D^{I}(t_1) \leq D^{I}(t_2) \leq D^{I}(t_3) \leq \ldots \leq D^{I} (t_T).  \]
Due to the permutation invariance of the problem and the monotonicity of the cost function we expect that:
\[   d^*(t_1) \geq d^*(t_2) \geq \ldots \geq d^*(t_T). \]
Too see this, assume by contradiction that for some $t_k$ we have
$d^*(t_k)< d^*(t_{k+1})$. Then, $D^{I}(t_k) + d^*(t_k) < D^{I}(t_{k+1}) +
d^*(t_{k+1})$. Hence reducing $d(t_{k+1})$ by a small amount of power $\Delta$ (equal for instance to $(d^*(t_{k+1})- d^*(t_k))/2$ and increasing $d(t_{k})$ by the same amount, yields a strict price reduction.
Let us denote $\tilde{d}$ this new demand profile, which fulfills:
\[ \tilde{d} (t) = \left \{  \begin{array}{rl}
d^*(t) & \textrm{if } t \notin \{ t_k, t_{k+1} \} \\
\frac{d^*(t_k)+d^*(t_{k+1})}{2} &
\textrm{if } t \in \{t_k, t_{k+1} \}
\end{array} \right .
 \]
Feasibility is also preserved since, any constraint involving both $t_k$ and $t_{k+1}$ is unaltered by the power swap. Any constraint
involving only $t_{k+1}$ is also fulfilled since $0 \leq \tilde{d}(t_{k+1}) \leq d^*(t_{k+1})$. Every constraint involving only $t_k$ is also fulfilled since the time $t_k$ can be swapped with time $t_{k+1}$:
\begin{equation*}
\begin{aligned}
&\sum_{t \in W} \tilde{d}(t) = \sum_{t \in W \cup \{t_{k+1} \} \backslash \{t_k \} } \tilde{d}(t) \\ &\leq
\sum_{t \in W \cup \{t_{k+1} \} \backslash \{t_k \} } d^*(t)
 \leq \sum_{j \in \mathcal{N}} \min \{  \textrm{card}(W), \bar{E}_j/ \bar{P}_j   \} \bar{P}_j.
\end{aligned}
\end{equation*}
% \[ \sum_{t \in W} \tilde{d}(t) = \sum_{t \in W \cup \{t_{k+1} \} \backslash \{t_k \} } \tilde{d}(t) \leq
% \sum_{t \in W \cup \{t_{k+1} \} \backslash \{t_k \} } d^*(t)
%  \leq \sum_{j \in \mathcal{N}} \min \{  \textrm{card}(W), \bar{E}_j/ \bar{P}_j   \} \bar{P}_j. \]
Hence, by convexity of the set $\mathcal{D}$, this contradicts optimality of $d^*(t)$ for all sufficiently small power swap among $t_k$ and $t_{k+1}$ in the direction of $\tilde{d}$ considered above. \\
Notice that, for any $W$ of cardinality $k$ we have the following:
\[   \sum_{t \in W} d^*(t)  \leq \sum_{i=1}^k d^*(t_i),   \] 
hence the same optimal generation cost and solution can be computed by solving:
\begin{subequations}\label{modifiedaggregate}
\begin{align}
\min_{\begin{aligned}
	d (t)\,\,&  t \in \mathcal{T} \\ g_i(t)\,\,& i \in \mathcal{I}, t \in \mathcal{T} 
	\end{aligned} } & \sum_{t \in \mathcal{T}} \sum_{i \in \mathcal{I} }   C_i (g_i(t)) \tag{\ref{modifiedaggregate}}\\
s.t.\quad \sum_{i \in \mathcal{I}} g_i(t) & = D^{I}(t) + d(t) \quad \forall \, t \in \mathcal{T}  \\ 
\underline{G}_i \leq & g_i(t) \leq \bar{G}_i \quad \forall \, t \in \mathcal{T}, \; \forall \, i \in \mathcal{I} \\
\sum_{t \in \mathcal{T}} d(t) & = \sum_{j \in \mathcal{N}} \bar{E}_j \\
\sum_{i=1}^{k} d(t_i) & \leq \sum_{j \in \mathcal{N}}
\min \{ k, \bar{E}_j / \bar{P}_j \} \bar{P}_j \nonumber \\ & \qquad \qquad\qquad\,\, \forall \, k \in \{1,2, \ldots, T\}. 
\end{align}
\end{subequations}

% \[
% \min_{\begin{array}{rl}
% 	d (t)&  t \in \mathcal{T} \\ g_i(t) & i \in \mathcal{I}, t \in \mathcal{T} 
% 	\end{array} }  \sum_{t \in \mathcal{T}} \sum_{i \in \mathcal{I} }   C_i (g_i(t)) \] 
% \[ \textrm{subject to} \]
% \[
% \begin{array}{rl}
% \sum_{i \in \mathcal{I}} g_i(t) &= D^{I}(t) + d(t) \qquad \forall \, t \in \mathcal{T}  \\ 
% \underline{G}_i \leq g_i(t) \leq & \bar{G}_i \qquad \forall \, t \in \mathcal{T}, \; \forall \, i \in \mathcal{I} \\
% \sum_{i=1}^{k} d(t_i) & \leq \sum_{j \in \mathcal{N}}
% \min \{ k, \bar{E}_j / \bar{P}_j \} \bar{P}_j \qquad \forall \, k \in \mathcal{T} \\
% \sum_{t \in \mathcal{T}} d(t) & = \sum_{j \in \mathcal{N}} \bar{E}_j   
% \end{array}   
% \]	
\section{Aggregate models for partial availability}\label{sec:partialavailability}
Next we aim to address the more general case of batteries with partial availability, viz. 
$\mathcal{A}_j \subseteq \mathcal{T}$ with strict inclusion also allowed.
This is of great practical relevance as it corresponds to batteries which are physically connected to the grid only during a limited amount of time, i.e. as typical of EVs. 

%Consider the following minimisation problem:
%\[   \mathcal{W}_{\min}  :=   \arg \min_{W \subset \mathcal{T} } \bar{D} (W).  \]

Unlike for the case of full availability, the expected ordering of total demand levels (and generation prices) cannot be inferred by looking at inflexible demand alone.
This is because different time intervals typically exhibit different sets of available batteries, and this in turn impacts ability and convenience to charge at any given time slot.
The next Lemma plays an important role to gain an a priori understanding of optimal power profiles.

\begin{lemma}
	\label{optimalagent}
	Let $u_j^*$ be optimal dispatch policies for problem (\ref{original})
	corresponding to the optimal aggregate demand profile $d^*(t)$ for problem (\ref{modified}), where
	\[  \sum_{j \in \mathcal{N}} u_j^*(t) = d^*(t) \qquad \forall \, t \in \mathcal{T}. \]
	Then, for all $j \in \mathcal{N}$, there exists $d_j >0$ such that:
	\begin{equation}
	\label{lemmaclaim}
	u_j^*(t) = \left \{ \begin{array}{rl} \bar{P}_j & \textrm{if } t \in \mathcal{A}_j \textrm{ and } D^{I}(t) + d^*(t) < d_j \\
	0 & \textrm{if } t \in \mathcal{A}_j \textrm{ and } D^{I}(t) + d^*(t) > d_j.   \end{array} \right .
	\end{equation}	
\end{lemma}

We say that $W \subset \mathcal{T}$ corresponds to an active constraint for $d^*(t)$ if:
\[  \sum_{t \in W} d^*(t) = \sum_{j \in \mathcal{N}} \min \{
\textrm{card} (\mathcal{A}_j \cap W), \bar{E}_j / \bar{P}_j   \} \bar{P}_j.   \]
Our aim is to devise an efficient algorithmic procedure to a priori select the support of active constraints in (\ref{modified}) so as to possibly retain the same optimal solution.
The following Lemma provides an important insight on some of the active constraints.
 
\begin{lemma}
	\label{activeconstraints}
	Let $d^*(t)$ be an optimal aggregated demand profile for problem (\ref{modified}). Consider for any $d > 0$, the following sets:
	\begin{equation}
	\label{thesublevelsets}
	  W_d := \{ t \in \mathcal{T}:  D^{I}(t) + d^* (t) < d  \}. 
	  \end{equation}
	Any such set is the support of an active constraint, viz.
	\[   \sum_{t \in W_d}  d^*(t) = \sum_{j \in \mathcal{N} }
	\min \{ \textrm{card} ( \mathcal{A}_j \cap W_d ), \bar{E}_j/ \bar{P}_j \} \bar{P}_j.  \] 
\end{lemma}

While Lemma \ref{activeconstraints} identifies an important set of active constraints, it does not provide any indication of how to compute these without solving the optimisation problem (\ref{modified}). In fact, the definition of $W_d$ explicitly requires the optimal aggregate power profile $d^*(t)$.
 The next result gives some important hints on how to look for active constraints by only using a priori available information. \\
To this end define, for each $W \subset \mathcal{T}$, the maximum average total demand as follows:
\begin{equation*}
\begin{aligned}
&\bar{D} (W) := \\   &\frac{ \sum_{t \in W} D^{I}(t)  + \sum_{t \in W} \sum_{j \in \mathcal{N}} \min \{
	\textrm{card} (\mathcal{A}_j \cap W), \bar{E}_j / \bar{P}_j   \} \bar{P}_j   }{\textrm{card}(W) }. 
\end{aligned}
\end{equation*}
% \[ \bar{D} (W) :=    \frac{ \sum_{t \in W} D^{I}(t)  + \sum_{t \in W} \sum_{j \in \mathcal{N}} \min \{
% 	\textrm{card} (\mathcal{A}_j \cap W), \bar{E}_j / \bar{P}_j   \} \bar{P}_j   }{\textrm{card}(W) }.   \]
Notice that this is affected both by the inflexible demand and by the availability windows of individual agents, and is computed assuming that within $W$ and compatibly with individual availability windows, batteries are drawing energy at their maximal rates (until completion of the charging task). 
 
In particular, Lemma \ref{achieved} states that $W_d$ achieves the minimum value of maximum average available demand $\bar{D}(W)$ over sets of its cardinality.
\begin{lemma}
	\label{achieved}
	Let $d^*(t)$ be the optimal aggregated demand for problem (\ref{modified}) and $d \geq 0$ be arbitrary. Define $W_d$
	according to (\ref{thesublevelsets}). Then the following holds:
\begin{equation}
\label{minimumcard}
    \bar{D} (W_d) =  \min_{W \subset \mathcal{T}: \textrm{card} (W)= \textrm{card}(W_d) }   \bar{D} (W).  
\end{equation}    	
Moreover $W_d$ is the unique minimizer in (\ref{minimumcard}).
\end{lemma}

While the $W_d$ sets defined in (\ref{thesublevelsets}) are not the only active constraints that one may find by solving (\ref{modified}), it turns out that they are enough to fully characterized the optimal solution $d^*(t)$, as clarified in the following Lemma \ref{wdenough}.

\begin{lemma}
	\label{wdenough}
	Let $d^*(t)$ be any optimal profile of problem (\ref{modified}) and 
	$W_d$ be defined as in (\ref{thesublevelsets}).
	Then the minimum of (\ref{modified}) is equal to the minimum of:
	\begin{subequations}\label{modified2}
	\begin{align}
	\min_{\begin{aligned}
		\tilde{d}(t)\,\,&  t \in \mathcal{T} \\ g_i(t)\,\,& i \in \mathcal{I}, t \in \mathcal{T} 
		\end{aligned} } &  \sum_{t \in \mathcal{T}} \sum_{i \in \mathcal{I} }   C_i (g_i(t)) \tag{\ref{modified2}}\\
	s.t. \quad  \sum_{i \in \mathcal{I}} g_i(t) &= D^{I}(t) + \tilde{d}(t) \quad \forall \, t \in \mathcal{T}  \\
	\underline{G}_i \leq & g_i(t) \leq \bar{G}_i \quad \forall \, t \in \mathcal{T}, \; \forall \, i \in \mathcal{I} \\
	\sum_{t \in \mathcal{T}} \tilde{d}(t) & = \sum_{j \in \mathcal{N}} \bar{E}_j \\
	\sum_{t \in W_d} \tilde{d}(t) & \leq \sum_{j \in \mathcal{N}}
	\min \{ \textrm{card} ( \mathcal{A}_j \cap W_d), \bar{E}_j / \bar{P}_j \} \bar{P}_j \nonumber \\ & \qquad \qquad \qquad \qquad\qquad \quad \forall \, d \geq 0 .
	\end{align}
	\end{subequations}
% 	\begin{equation}
% 	\label{modified2}
% 	\min_{\begin{array}{rl}
% 		\tilde{d} (t)&  t \in \mathcal{T} \\ g_i(t) & i \in \mathcal{I}, t \in \mathcal{T} 
% 		\end{array} }  \sum_{t \in \mathcal{T}} \sum_{i \in \mathcal{I} }   C_i (g_i(t)) 
% 	\end{equation} 
% 	\[ \textrm{subject to} \]
% 	\[
% 	\begin{array}{rl}
% 	\sum_{i \in \mathcal{I}} g_i(t) &= D^{I}(t) + \tilde{d}(t) \qquad \forall \, t \in \mathcal{T}  \\ 
% 	\underline{G}_i \leq g_i(t) \leq & \bar{G}_i \qquad \forall \, t \in \mathcal{T}, \; \forall \, i \in \mathcal{I} \\
% 	\sum_{t \in W_d} \tilde{d}(t) & \leq \sum_{j \in \mathcal{N}}
% 	\min \{ \textrm{card} ( \mathcal{A}_j \cap W_d), \bar{E}_j / \bar{P}_j \} \bar{P}_j \qquad \forall \, d \geq 0  \\
% 	\sum_{t \in \mathcal{T}} \tilde{d}(t) & = \sum_{j \in \mathcal{N}} \bar{E}_j 
% 	\end{array}   
% 	\]	
\end{lemma}
It is worth pointing out that as $d$ ranges in $[0,+\infty)$ there are at most $T$ distinct values of $W_d$. This is because the $W_d$ are nested with each other. We may therefore rename them and denote in ascending order as:
\[ W_1 \subsetneq W_2 \subsetneq W_3 \subsetneq \ldots \subsetneq W_Q \]
for some positive integer $Q \leq T$.

Let $d^*(t)$ be any optimal profile of problem (\ref{modified}) and 
$W_d$ be defined as in (\ref{thesublevelsets}).
We already remarked that there are at most $T$ distinct values of $W_d$ and that they are nested with one another so that we may rename them and denote in ascending order as:
\[ W_1 \subsetneq W_2 \subsetneq W_3 \subsetneq \ldots \subsetneq W_Q \]
Moreover, it is clear by their definition that the following implication holds:
\[     \frac{\sum_{t \in W_k}  D^{I}(t) + d^*(t)}{\textrm{card} (W_k)}  < \frac{ \sum_{t \in W_{k+1}} D^{I}(t) + d^*(t)  }{\textrm{card}(W_{k+1})}.  \]
As a consequence:
\[    \bar{D} (W_k) < \bar{D} ( W_{k+1}). \]
The discussion so far clearly shows that
the following procedure can be adopted to pre-compute $T$ constraints to be included in a simplified version of (\ref{modified}).
Define $W^*_k$ as follows:
\begin{equation}
\label{cardinality}
W^*_k := \arg \min_{W \subset \mathcal{T}: \textrm{card}(W) = k } \bar{D} (W). 
\end{equation}
Notice that the sets $W^*_k$, unlike the $W_d$ previously discussed, can be computed on the basis of a priori available information. Moreover,  by virtue of
Lemma \ref{achieved}, the following inclusion holds:
\begin{equation}
\label{nested}
     \{ W_d: d \geq 0 \} \subset \{ W^*_k: k \in \{1,2,\ldots,T \} \} \subset 2^{\mathcal{T}}. 
\end{equation}     
Hence, we may replace the constraints of support $W_d$ in (\ref{modified2}) with similar constraints of support 
$W^*_k$, without affecting the cost of the optimal solution. 
This is stated below and is our main result for this section: 
\begin{theorem}
		\label{wstarenough}
		Let $d^*(t)$ be any optimal profile of problem (\ref{modified}) and 
		$W^*_k$ be defined as in (\ref{cardinality}).
		Then the minimum of (\ref{modified}) is equal to the minimum of:
		\begin{subequations}\label{modified3}
		\begin{align}
		\min_{\begin{aligned}
			\tilde{d}(t)\,\,&  t \in \mathcal{T} \\ g_i(t)\,\,& i \in \mathcal{I}, t \in \mathcal{T} 
			\end{aligned} } & \sum_{t \in \mathcal{T}} \sum_{i \in \mathcal{I} }   C_i (g_i(t)) \tag{\ref{modified3}}\\
			s.t. \quad \sum_{i \in \mathcal{I}} g_i(t) &= D^{I}(t) + \tilde{d}(t) \quad \forall \, t \in \mathcal{T}  \\ 
		\underline{G}_i \leq & g_i(t) \leq \bar{G}_i \quad \forall \, t \in \mathcal{T}, \; \forall \, i \in \mathcal{I} \\
		\sum_{t \in \mathcal{T}} \tilde{d}(t) & = \sum_{j \in \mathcal{N}} \bar{E}_j \\
		\sum_{t \in W^*_k} \tilde{d}(t) & \leq \sum_{j \in \mathcal{N}}
		\min \{ \textrm{card} ( \mathcal{A}_j \cap W^*_k), \bar{E}_j / \bar{P}_j \} \bar{P}_j \nonumber \\ & \qquad \qquad\qquad\,\, \forall \, k \in \{1,2, \ldots, T\}.   
		\end{align}
		\end{subequations}
% 		\begin{equation}
% 		\label{modified3}
% 		\min_{\begin{array}{rl}
% 			\tilde{d} (t)&  t \in \mathcal{T} \\ g_i(t) & i \in \mathcal{I}, t \in \mathcal{T} 
% 			\end{array} }  \sum_{t \in \mathcal{T}} \sum_{i \in \mathcal{I} }   C_i (g_i(t)) 
% 		\end{equation} 
% 		\[ \textrm{subject to} \]
% 		\[
% 		\begin{array}{rl}
% 		\sum_{i \in \mathcal{I}} g_i(t) &= D^{I}(t) + \tilde{d}(t) \qquad \forall \, t \in \mathcal{T}  \\ 
% 		\underline{G}_i \leq g_i(t) \leq & \bar{G}_i \qquad \forall \, t \in \mathcal{T}, \; \forall \, i \in \mathcal{I} \\
% 		\sum_{t \in W^*_k} \tilde{d}(t) & \leq \sum_{j \in \mathcal{N}}
% 		\min \{ \textrm{card} ( \mathcal{A}_j \cap W^*_k), \bar{E}_j / \bar{P}_j \} \bar{P}_j \qquad \forall \, k \in \{1,2, \ldots, T\}   \\
% 		\sum_{t \in \mathcal{T}} \tilde{d}(t) & = \sum_{j \in \mathcal{N}} \bar{E}_j 
% 		\end{array}   
% 		\]	
\end{theorem} 

\begin{algorithm}
\DontPrintSemicolon 

\textbf{Initialization:}

\nonl $k = T$

\textbf{Update procedure:}

\While{$k > 0$}
{

Find $W^{*}_{k}$ by solving \eqref{cardinality}

Select the inequality corresponding to $W^{*}_{k}$ as an active constraint, decrease the cardinality by
$k = k - 1$
}

\textbf{Final results:} The set of selected active constraints are returned.
\caption{Combinatorial Optimisation Algorithm}
\label{CombinatorialAlgorithm}
\end{algorithm}

The active constraint selection is described in Algorithm \ref{CombinatorialAlgorithm}. One drawback of (\ref{cardinality}) is that, for every assigned value of $k$, it is a combinatorial optimisation problem that entails generation of $\binom{T}{k}$
%%\left ( \begin{array}{c} T \\ k \end{array} \right)$ 
sets and computation of the corresponding maximum available power.
 While this would typically be faster than a linear programming solver dealing with $2^T$ constraints, it might still be computationally inefficient, and lead to scalability issues.
Hence, we propose in the following two greedy algorithms (of limited computational complexity) for the pre-computation of active constraints, which in simulation have achieved near-optimal performance for large classes of randomly generated problems.
 
%The monotonicity of $\bar{D}(W_k)$ with respect to index $k$ can be exploited to propose an algorithm for a priori computation of the set of active constraints.

\begin{algorithm}
\DontPrintSemicolon 

\textbf{Initialization:}

\nonl $k = T; \hat{W}_T = \mathcal{T}; S = 1$

\textbf{Update procedure:}

\While{$k > 1$}
{

Calculate $\hat{W}_{k-S} := \arg \min_{W \subset \hat{W}_k, |W|=k-S} \bar{D} (W)$

\eIf{ $\hat{W}_{k-S}$ is not a singleton or $\bar{D} (\hat{W}_{k-S}) \geq \bar{D} (\hat{W}_k)$}{$S = S + 1$
}
{

Select the inequality corresponding to $\hat{W}_{k-S}$ as an active constraint, decrease the cardinality by $k = k - S$, and reset $S = 1$

}

}

\textbf{Final results:} The set of selected $\leq T$ active constraints are returned.
\caption{Greedy Algorithm I}
\label{GreedyAlgorithm1}
\end{algorithm}

% \begin{tabular}{|l|}
% \hline
% \\
% \textbf{Greedy Algorithm 1} \\
% \\
% 1. Initialization: $k:=T$; $\hat{W}_T:= \mathcal{T}$; \\
% 2. $S=1$ \\
% 3. $\hat{W}_{k-S} := \arg \min_{W \subset \hat{W}_k, |W|=k-S} \bar{D} (W)$; \\
% 4. If $\hat{W}_{k-S}$ is not a singleton or $\bar{D} (\hat{W}_{k-S}) \geq \bar{D} (\hat{W}_k)$ then increase $S$ and go to $3$. \\
% 5. Decrease $k$ by $S$; if $k>0$ then go to 2. \\
% \\
% \hline
% \end{tabular} \\

The greedy algorithm I (in Algorithm \ref{GreedyAlgorithm1}) constructs sets $\hat{W}_k$ in descending order, nested within one another, and starting from the set of cardinality $T$.
It avoids a combinatorial explosion as it seeks to minimize $\bar{D}(W)$ only generating sets which are nested in the previously determined minimal set.
However, it is not guaranteed to determine all of the $W^*_k$ sets and consequently it is also not guaranteed to generate all of the $W_d$ ($d \geq 0$) sets.

In most cases, the greedy algorithm I searches for the minimal $\bar{D}(W)$ over a relatively small number of candidate solutions, i.e., $\binom{k}{k-S}$, when $S$ is small. However, large values of $S$ could slow down the algorithm, due to the unavoidable combinatorial explosion. To further release the computational burden, we propose an alternative greedy algorithm II in Algorithm \ref{GreedyAlgorithm2} for active constraints selection as follows:

\begin{algorithm}
\DontPrintSemicolon 

\textbf{Initialization:}

\nonl $k:=T; \hat{W}_T:= \mathcal{T}$

\textbf{Update procedure:}

\While{$k > 1$}
{

Calculate $\hat{W}^{(m)}_{k-1} := \arg \min_{W \subset \hat{W}^{(m)}_k, |W|=k-1} \bar{D} (W)$, $\forall \hat{W}^{(m)}_k \in \hat{W}_k = \{\hat{W}^{(1)}_k, \hat{W}^{(2)}_k, \cdots \hat{W}^{(M)}_k\}$

Let $\hat{W}_{k-1} = \arg \min_{\hat{W}^{(m)}_{k-1} \in \{\hat{W}^{(1)}_{k-1}, \hat{W}^{(2)}_{k-1}, \cdots, \hat{W}^{(M)}_{k-1}\}} \bar{D}(\hat{W}^{(m)}_{k-1})$

Select the inequality corresponding to $\hat{W}_{k-1}$ as an active constraint

Decrease the cardinality by $k = k - 1$

}

\textbf{Final results:} The set of selected active constraints are returned.
\caption{Greedy algorithm II}
\label{GreedyAlgorithm2}
\end{algorithm}

% \begin{tabular}{|l|}
% \hline
% \\
% \textbf{Greedy Algorithm 2} \\
% \\
% 1. Initialization: $k:=T$; $\hat{W}_T:= \mathcal{T}$;\\
% 2. Calculate $\hat{W}^{(m)}_{k-1} := \arg \min_{W \subset \hat{W}^{(m)}_k, |W|=k-1} \bar{D} (W)$, $\forall \hat{W}^{(m)}_k \in \hat{W}_k = \{\hat{W}^{(1)}_k, \hat{W}^{(2)}_k, \cdots \hat{W}^{(M)}_k\}$;\\
% 3. Let $\hat{W}_{k-1} = \arg \min_{\hat{W}^{(m)}_{k-1} \in \{\hat{W}^{(1)}_{k-1}, \hat{W}^{(2)}_{k-1}, \cdots, \hat{W}^{(M)}_{k-1}\}} \bar{D}(\hat{W}^{(m)}_{k-1})$;\\
% 4. Decrease $k$ by $1$; if $k>0$ then go to 2.\\
% \\
% \hline
% \end{tabular}\\

The advantage of this algorithm is that the number of candidates $W$ at the $k$-th iteration is $\binom{k}{k-1} = k$ which avoids the combinatorial explosion and terminates in a number of steps of order $T^2$. However, this fast greedy algorithm might lose optimality for certain examples because the minimizer $\hat{W}_{k-1}$ is restricted to be selected from the subset of $\hat{W}_{k}$ which disables the active constraint selection for $S>1$ in greedy algorithm I.

\section{Advantages in computation}
\label{sec:computation}

This section is devoted to apply the proposed aggregate demand model based optimisation problem and the active constraints selection algorithms to some numerical case studies. All the tested optimisation problems were solved by the YALMIP toolbox and Gurobi optimiser in MATLAB R2019a, on a computer with 2-core 3.50GHz Intel(R) Xeon(R) E5-1650 processor and 32GB RAM.

% \subsection{Computational performance of proposed algorithms}
The computational advantages of using optimisation problem \eqref{modified} and the greedy constraint selection algorithms in a large population of storage devices are evaluated. In particular, the following methods for exploring the optimal dispatch of storage fleets are considered throughout the rest of this section
\begin{itemize}
    \item M1: Each storage unit is modeled as an individual agent and the optimisation problem \eqref{original} is solved;
    \item M2: Aggregate demand model based optimisation by solving problem \eqref{modified} with $2^{T}-1$ inequality constraints on aggregate demand;
    \item M3: Aggregate demand model based optimisation by solving problem \eqref{modified3} with $T$ inequality constraints defined in \eqref{cardinality};
    \item M4: Aggregate demand model based optimisation by solving problem \eqref{modified3} with $\leq T$ constraints chosen by greedy algorithm I.
    \item M5: Aggregate demand model based optimisation by solving problem \eqref{modified3} with active constraints chosen by greedy algorithm II.
\end{itemize}

\begin{figure}[!ht]
\centerline{\includegraphics[width=0.48\textwidth]{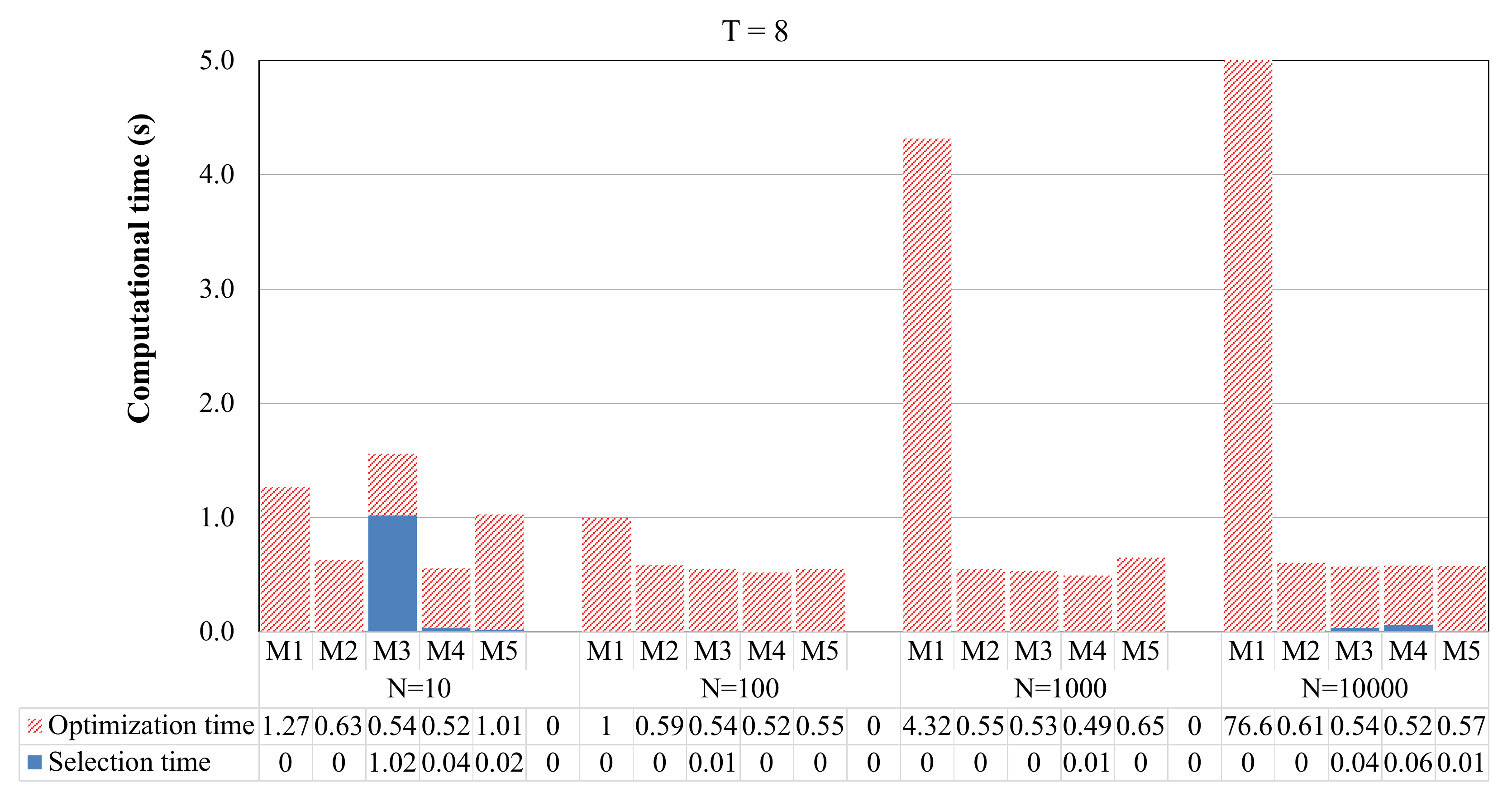}}
\vspace{-3mm}
\caption{Computational time of $T_{1} = 8$ h for different number of batteries.}
\label{DifferentN_FixedT8}
\end{figure}

\begin{figure}[!ht]
\centerline{\includegraphics[width=0.48\textwidth]{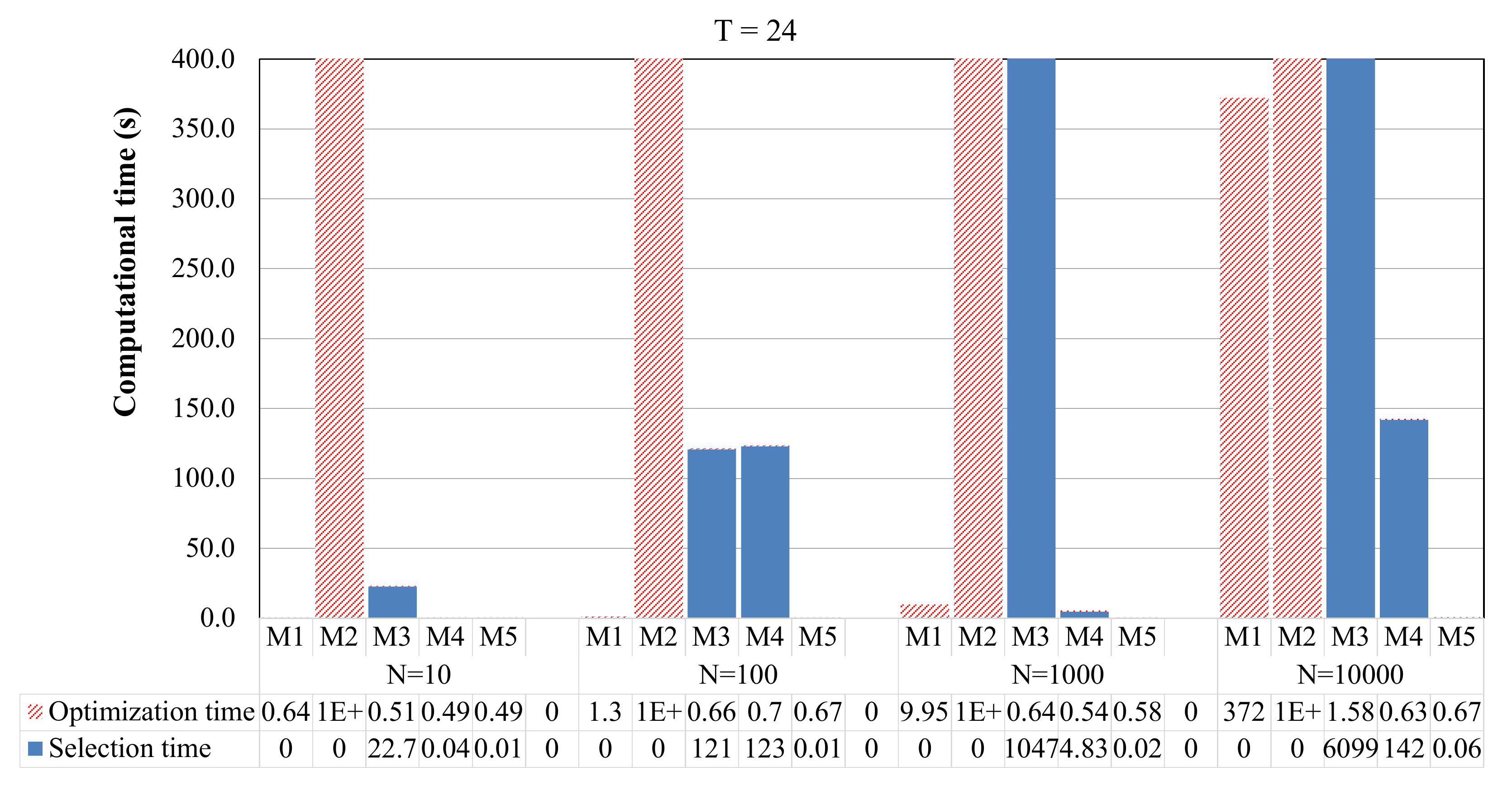}}
\vspace{-3mm}
\caption{Computational time of $T_{2} = 24$ h for different number of batteries.}
\label{DifferentN_FixedT24}
\end{figure}

To compare the performance of these methods, a series of case studies are designed considering various length of full time window and different number of storage devices in a fleet, i.e.,
\begin{equation*}
    T_{1} = 8\, \text{h},\,\, T_{2} = 16\, \text{h},\,\, T_{3} = 24\, \text{h}
\end{equation*}
\begin{equation*}
    N_{1} = 10,\,\, N_{2} = 100,\,\, N_{3} = 1000,\,\, N_{4} = 10 000.
\end{equation*}
The data are sampled every hour. Then, given a time window $\mathcal{T}_{i} = \{1, 2, \cdots, T_{i}\}$ where $i \in \{1, 2, 3\}$, the availability window of each device is uniformly generated as a subset of $\mathcal{T}_{i}$. For each storage, the rated power is set to be $\bar{P}_{j} = 1$ kW  and the target energy $E_{j}$ follows a uniform distribution on the continuous interval $[0,  \bar{P}_{j} \cdot \textrm{card}(\mathcal{A}_{j})]$. Based on these randomly generated examples, the following results are observed:

The computational performances of M1 - M5 for $T_{1} = 8$ h and $T_{3} = 24$ h are shown in Fig.~\ref{DifferentN_FixedT8} and Fig.~\ref{DifferentN_FixedT24}. From the figures, it emerges that the time for solving optimisation problem \eqref{original}  is proportional to the number of agents, so M1 is inefficient for problems with large number of agents. In Fig.~\ref{DifferentN_FixedT8}, when the considered time window is of small length, approaches M2-M4 which are based on the aggregate demand model are effective regardless of the number of devices. In particular, these methods are much faster than M1 for $N_{4} = 10000$. However, the aggregate demand based methods may be exposed to high computational burden if the window is large as shown in Fig.~\ref{DifferentN_FixedT24}, i.e., $T_{3} = 24$. Specifically, M2 can provide the optimal aggregate demand profile, but it is computationally expensive due to $2^{24}-1 \approx 1.6 \cdot 10^7$ constraints. Since the active constraint selection in M3 - M5 is determined by the maximum average demand, a larger number of batteries can also increase the time spent on constraint selection, however, the solution of optimisation \eqref{modified3} is extremely fast and unaffected by the agents' number. Moreover, as described in Section \ref{sec:partialavailability}, combinatorial explosion increases computational burden in searching for the optimal $W^{*}_{k}$ in M3 and possibly $\hat{W}_{k-S}$ in M4, the time on active constraint selection is much higher than that of using the greedy algorithm II in M5.

\begin{figure}[!ht]
\centerline{\includegraphics[width=0.48\textwidth]{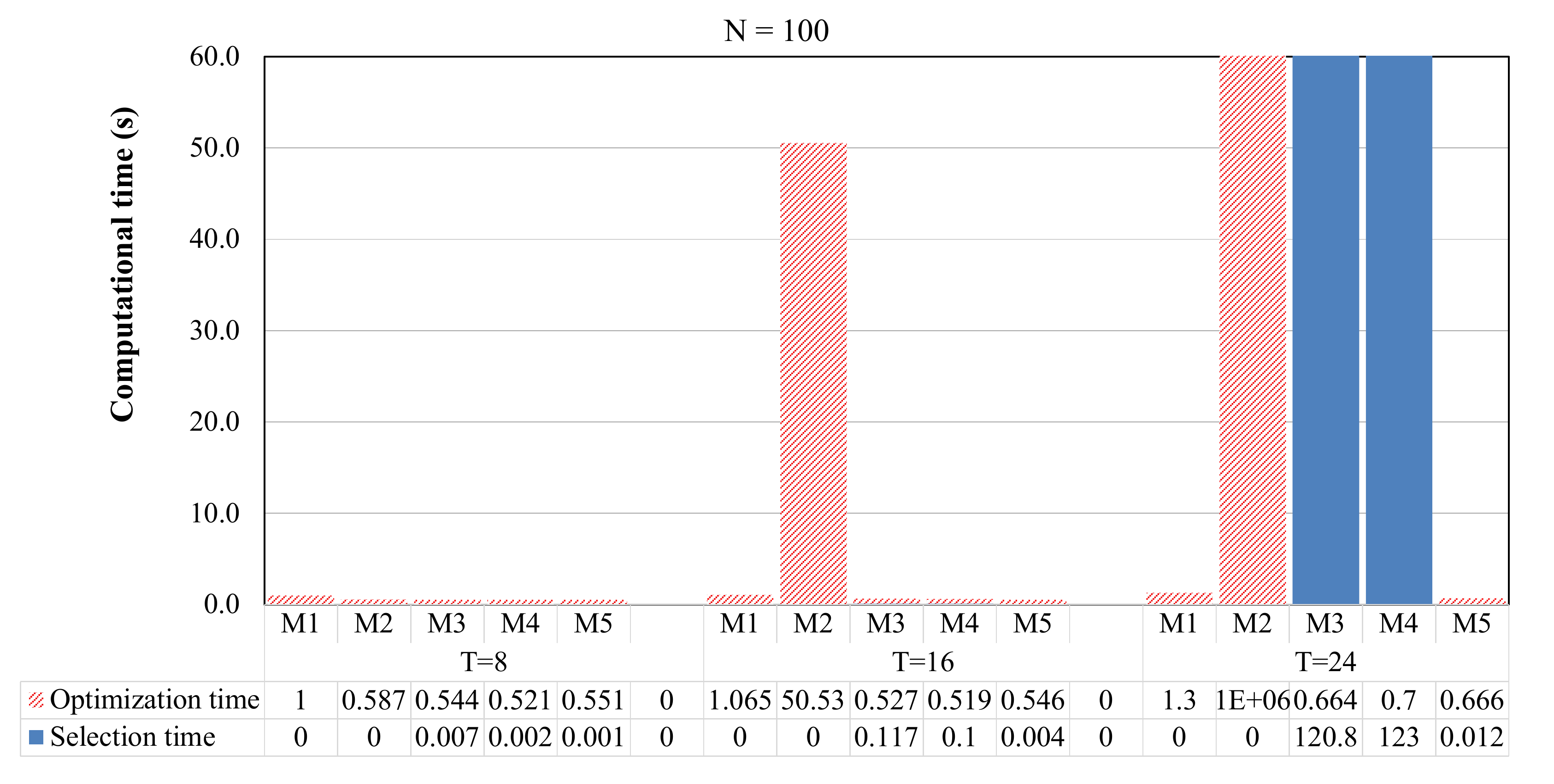}}
\vspace{-3mm}
\caption{Computational time of $N_{1} = 100$ h for different time horizon.}
\label{DifferentT_FixedN100}
\end{figure}

\begin{figure}[!ht]
\centerline{\includegraphics[width=0.48\textwidth]{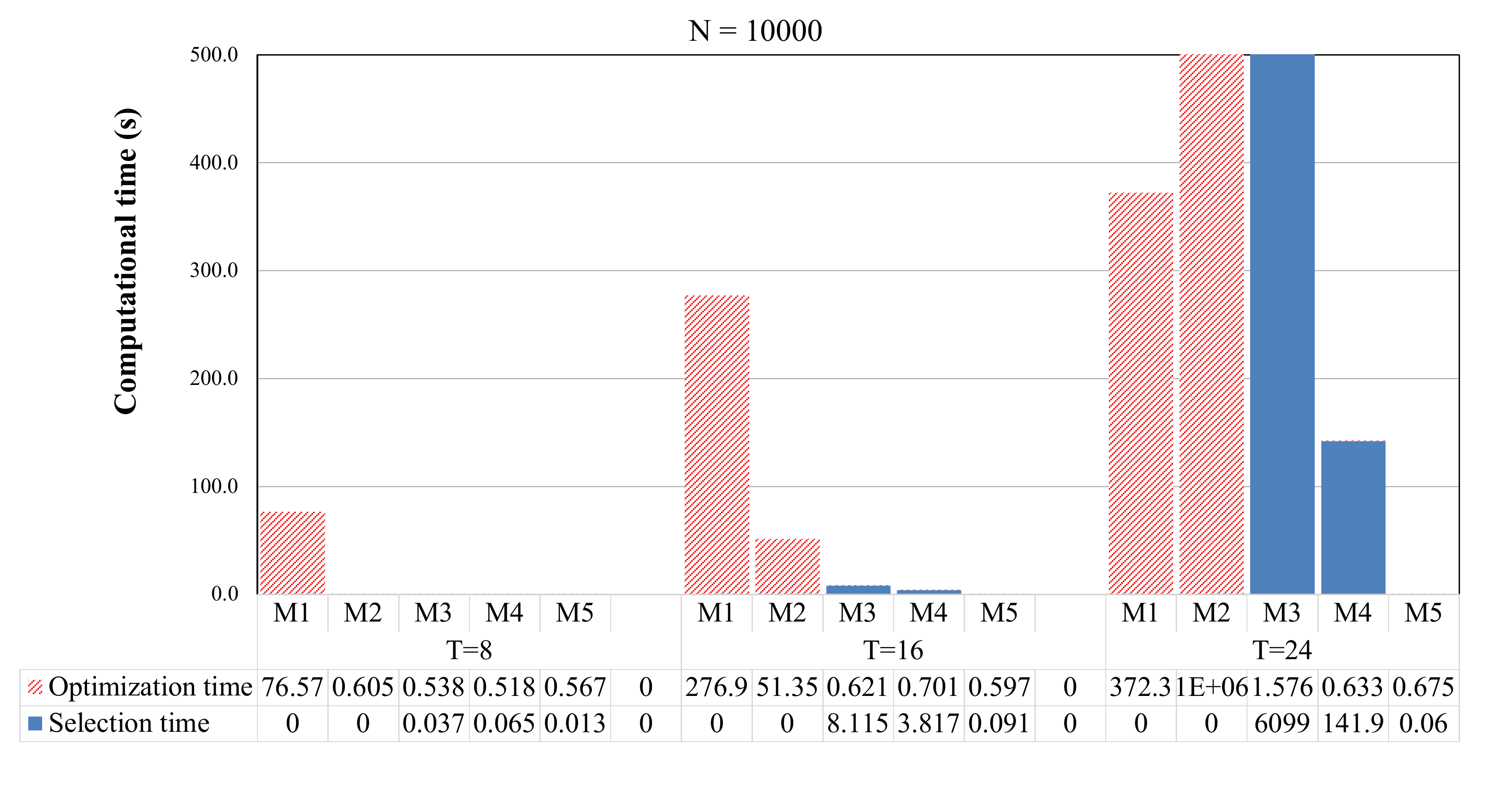}}
\vspace{-3mm}
\caption{Computational time of $N_{1} = 10000$ h for different time horizon.}
\label{DifferentT_FixedN10000}
\end{figure}

Fig.~\ref{DifferentT_FixedN100} and Fig.~\ref{DifferentT_FixedN10000} presents the computational time of each optimal power dispatch algorithm for $T_{1} - T_{3}$ with $N_{2}$ and $N_{4}$ batteries. The results of M1 verify that the computational complexity for solving optimisation problem \eqref{original} is linear with respect to $N \cdot T$. For the methods M3 and M4 which are based on aggregate demand model, the performance is mainly affected by the length of horizon and example parameters. Particularly, comparing the selection time of M4 at $T_{3} = 24$ h in both figures, we can conclude that less combinatorial explosion is encountered for $N_{4} = 10000$.

Notice that the computational performance of M5 using greedy algorithm II for active constraint selection is the best for problems with large number of batteries and long horizon interval. However, this improvement in computation is achieved at the expense of optimality for some examples. Although the methods M1 - M3 are suffering from the computational challenge, the results are the optimal. To quantify the performance in obtaining the optimal solution using M4 and M5, $10^{4}$ scenarios with $N_{1} = 10$ batteries and $T_{4} = 24$ h time horizon are randomly generated. Finally, the success rate to reach the optimal solution for M4 and M5 are 95.6\% and 93.6\%, respectively (with small optimality gaps in the case of missed optimal solutions).

\section{Two-area system}\label{sec:twoarea}

\begin{figure}[H]
\vspace{-3mm}
\centerline{\includegraphics[width=0.48\textwidth]{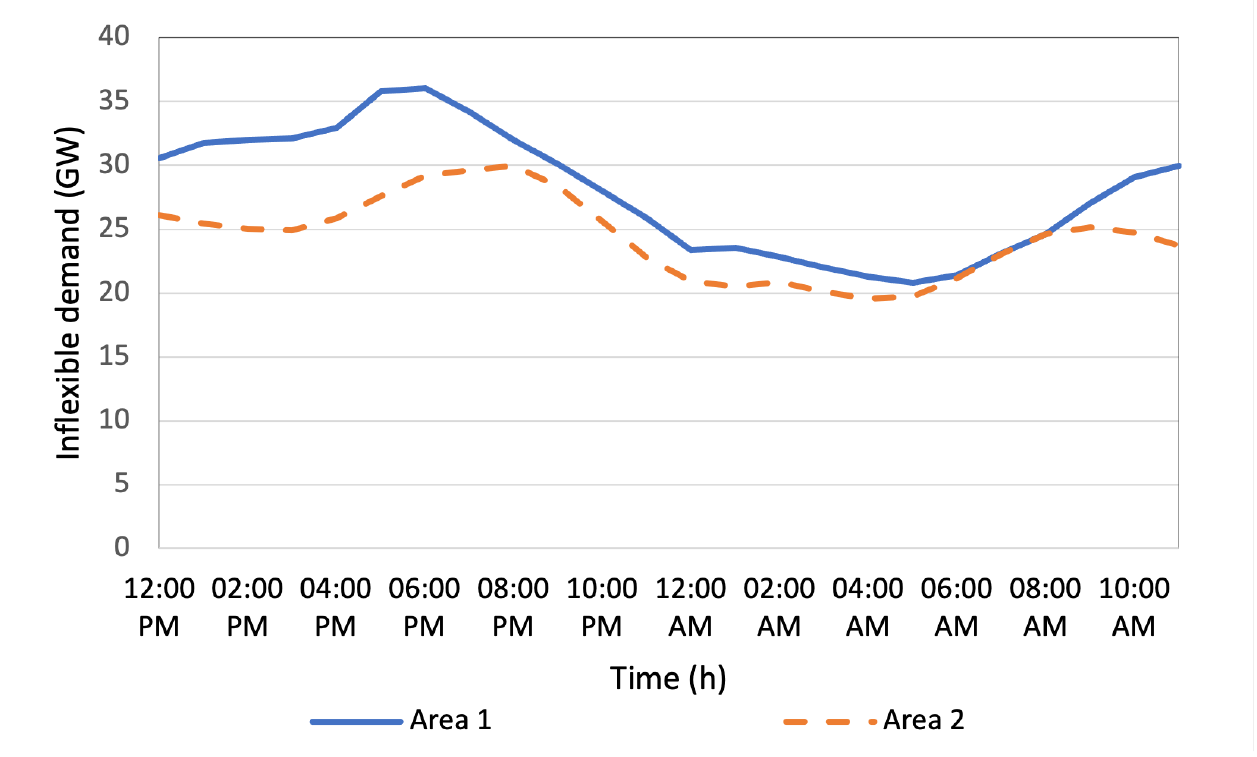}}
\vspace{-3mm}
\caption{Inflexible demand profiles at two areas.}
\label{InflexibleDemand}
\end{figure}

Next, we consider the application of the aggregate demand model on a two-area system $\mathcal{I} = \{1, 2\}$ over a time horizon $\mathcal{T} = \{1, 2, \cdots, 24\}$ h and an hourly discretization step. Each area is equipped with local dispatchable generators and numerous flexible storage devices. In this case study, to simplify the analysis, we assume that the size of a single generator is significantly smaller than the total installed power output capacity, i.e., the dispatchable power generation at area $i \in \mathcal{I}$ can be approximated by one generator whose power is a continuous variable between $\underline{G}_{i} = 0$ GW and $\bar{G}_{i} = 60$ GW. The generation cost functions are expressed as
\begin{equation*}
    C_{i}(g_{i}(t)) = a_{i}\cdot (g_{i}(t))^{2} + b_{i}\cdot g_{i}(t), \quad \forall i \in \mathcal{I}
\end{equation*}
where 
\begin{equation*}
\begin{aligned}
    a_{1} = 1\times & 10^{4}\; \text{\pounds/GW}^{2}\text{h}, \; b_{1} = 1.5\times 10^{4}\; \text{\pounds/GWh},\\
    a_{2} = 2\times & 10^{4}\; \text{\pounds/GW}^{2}\text{h}, \; b_{2} = 1.4\times 10^{4}\; \text{\pounds/GWh}.
\end{aligned}
\end{equation*}
Moreover, the inflexible demands of these two areas $D^{I}_{1}$ (Area 1) and $D^{I}_{2}$ (Area 2) are assumed to be equal to historic data of total power consumption in the U.K. National Grid \cite{websitedata}, as shown in Fig.~\ref{InflexibleDemand}. 

A significant penetration of flexible demand is required for charging while minimizing the total generation cost. Specifically, there are $4 \times 10^{6}$ batteries in area 1 and $6 \times 10^{6}$ batteries in area 2, i.e., $\mathcal{N}_{1} = \{1, 2, \cdots, 4 \times 10^{6}\}, \mathcal{N}_{2} = \{1, 2, \cdots, 6 \times 10^{6}\}$. It is assumed that the rated power is equal for the whole population, with $\bar{P}_{j} = 5$ kW, $\forall j \in \mathcal{N}_{1} \cup \mathcal{N}_{2}$. Each device can charge only within a continuous time interval $[t^{s}_{j}, t^{s}_{j} + t^{d}_{j}]$ h, where the starting time $t^{s}_{j}$ and available time duration $t^{d}_{j}$ follow normal distributions, with the following mean and standard deviation:
\begin{equation*}
    \mu_{t^{s}} = 18:00 \,\,\text{h}, \,\, \sigma_{t^{s}} = 1 \,\,\text{h}, \,\, \mu_{t^{d}} = 10 \,\,\text{h}, \,\, \sigma_{t^{d}} = 2 \,\,\text{h}.
\end{equation*}
The energy required $\bar{E}_{j}$ for task completion is uniformly generated on the interval $[0, \bar{P}_{j}\cdot t^{d}]$. Moreover, the power transmission between these two area is bounded by the transmission capacity, i.e., $\abs{p(t)} \leq \bar{p}$.

With the parameters in the two-area system, the optimal power management problem to be solved is formulated as
\begin{subequations}
\begin{align*}
\min_{\begin{aligned}
	u_j(t),\;&g_i(t),\;p(t)  \\ j \in \mathcal{N}_{i},\,\, & i \in \mathcal{I},\,\, t \in \mathcal{T} 
\end{aligned} } &  \sum_{t \in \mathcal{T}} \left[  C_{1} (g_{1}(t)) + C_{2} (g_{2}(t)) \right]\\
s.t.\quad  \underline{G}_i \leq & g_i(t) \leq \bar{G}_i \quad \forall \, t \in \mathcal{T}, \; \forall \, i \in \mathcal{I} \\
0 \leq u_j(t) \leq \bar{P}_j,\,\, & \bar{E}_j = \sum_{t \in \mathcal{T}} u_j (t) \quad \forall \, t \in \mathcal{T}, \; \forall \, j \in \mathcal{N}_{1} \cup \mathcal{N}_{2}\\
u_j(t) &= 0 \quad \forall \, t \notin \mathcal{A}_j,\; \forall \, j \in \mathcal{N}_{1} \cup \mathcal{N}_{2}\\
g_{1}(t) + p(t) &= D^{I}_{1}(t) + \sum_{j \in \mathcal{N}_{1}} u_j(t) \quad \forall \, t \in \mathcal{T} \\
g_{2}(t) - p(t) &= D^{I}_{2}(t) + \sum_{j \in \mathcal{N}_{2}} u_j(t) \quad \forall \, t \in \mathcal{T}\\
-\bar{p} \leq & p(t) \leq \bar{p}
\end{align*}
\end{subequations}
which has an equivalent form as following using the aggregate demand model, viz.,
\begin{subequations}
\begin{align*}
\min_{\begin{aligned}
	d_{i}(t),\; g_i(&t),\; p(t)  \\ i \in \mathcal{I},\,\,& t \in \mathcal{T} 
	\end{aligned} }& \,\, \sum_{t \in \mathcal{T}} \left[   C_{1} (g_{1}(t)) + C_{2} (g_{2}(t)) \right]\\
s.t.\quad \underline{G}_i \leq & g_i(t) \leq \bar{G}_i \quad \forall \, t \in \mathcal{T}, \; \forall \, i \in \mathcal{I} \\
g_{1}(t) + p(t) &= D^{I}_{1}(t) + d_{1}(t) \quad \forall \, t \in \mathcal{T} \\
g_{2}(t) - p(t) &= D^{I}_{2}(t) + d_{2}(t) \quad \forall \, t \in \mathcal{T}\\
-\bar{p} \leq & p(t) \leq \bar{p}\\
\sum_{t \in \mathcal{T}} d_{i}(t) & = \sum_{j \in \mathcal{N}_{i}} \bar{E}_j \quad \forall i \in \mathcal{I} \\
\sum_{t \in \mathcal{W}} d_{i}(t) & \leq \sum_{j \in \mathcal{N}_{i}}
\min \{ \textrm{card} ( \mathcal{A}_j \cap \mathcal{W}), \bar{E}_j / \bar{P}_j \} \bar{P}_j \nonumber \\
& \qquad \qquad \qquad \quad  \quad  \forall \, \mathcal{W} \subset \mathcal{T}, \; \forall \, i \in \mathcal{I}. 
\end{align*}
\end{subequations}

As illustrated in Section \ref{sec:computation}, solving the above two optimisation problems directly is computational inefficient when the number of storage units is large or the time horizon is long. Thus, we adopt two heuristic constraint selection methods for the aggregate demand model to ease the computation burden. Specifically, all the considered approaches to minimize the total generation cost for the two-area power system are summarized as below, each of which is represented by an abbreviation for notational convenience,
\begin{itemize}
    \item MM1: Each storage unit is modeled as an individual agent whose operational behavior is the decision variable.
    \item MM2: The storage units in each area are aggregated as an aggregate demand model, and the power subset of full time windows are imposed in constraints formulation for both areas ($2\times(2^{T}-1)$ constraints).
    \item MM3: The storage units in each area are aggregated as an aggregate demand model, and $2T$ constraints are selected where each area selects its $T$ constraints according to the greedy algorithm II.
    \item MM4: The storage units in each area are aggregated as an aggregate demand model, and the $T$ constraints are selected by applying greedy algorithm II to the merged fleet $\mathcal{N}_1 \cup \mathcal{N}_2$.
\end{itemize}
Particularly, in methods MM3 and MM4 with reduced number of constraints, for the aim of feasibility satisfaction, $d_{i}(t) \leq \bar{d}_{i}(t) = \sum_{j \in \mathcal{N}_{i}: t \in \mathcal{A}_{j}} \bar{P}_{j}$ is also enforced, that is the aggregate demand $d_{i}(t)$ of both areas cannot exceed their maximum power rating.

\begin{figure}[H]
\centerline{\includegraphics[width=0.48\textwidth]{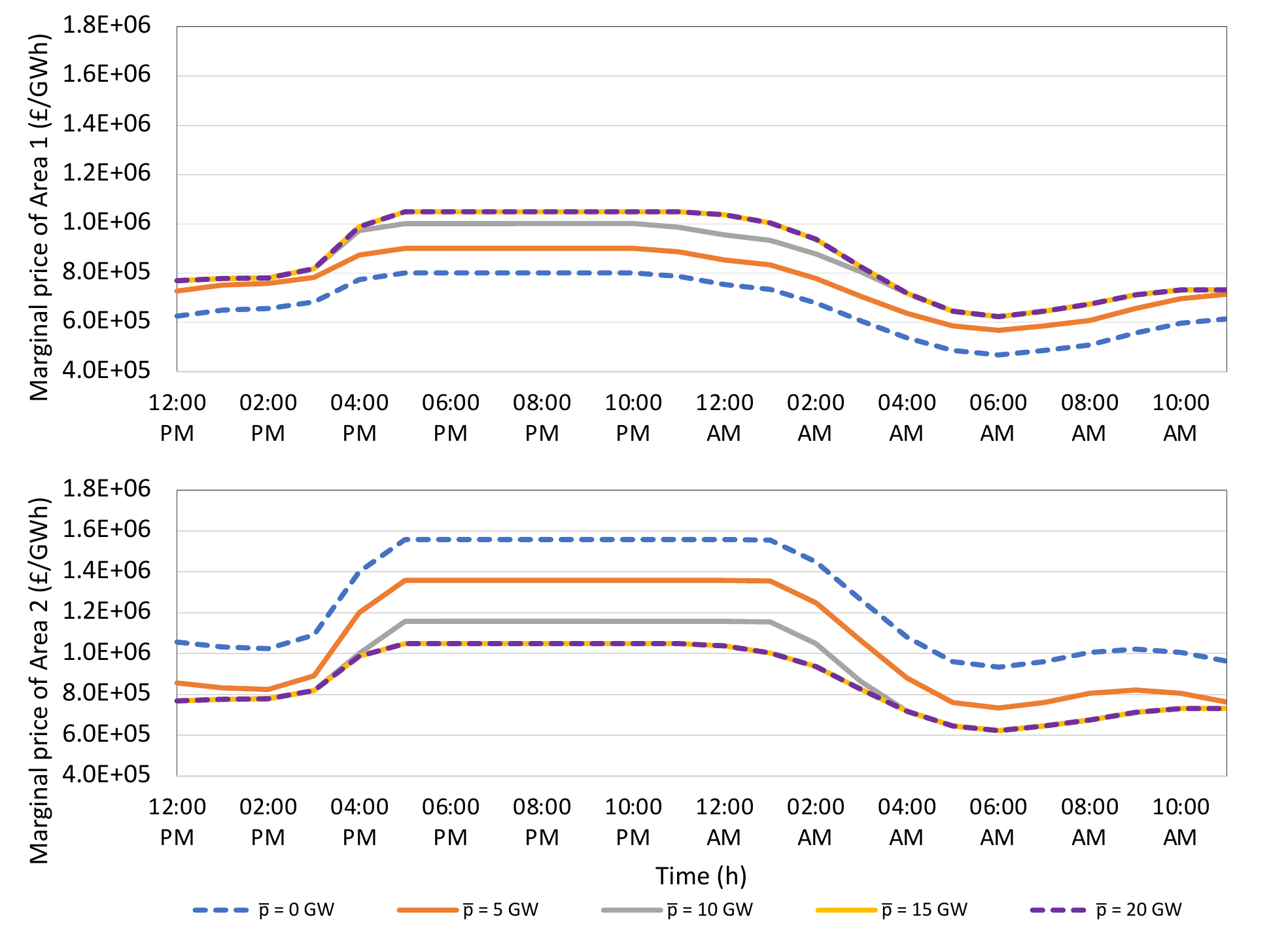}}
\vspace{-3mm}
\caption{The marginal price profiles at two areas.}
\label{MarginalPrices}
\end{figure}

\begin{figure}[H]
\centerline{\includegraphics[width=0.48\textwidth]{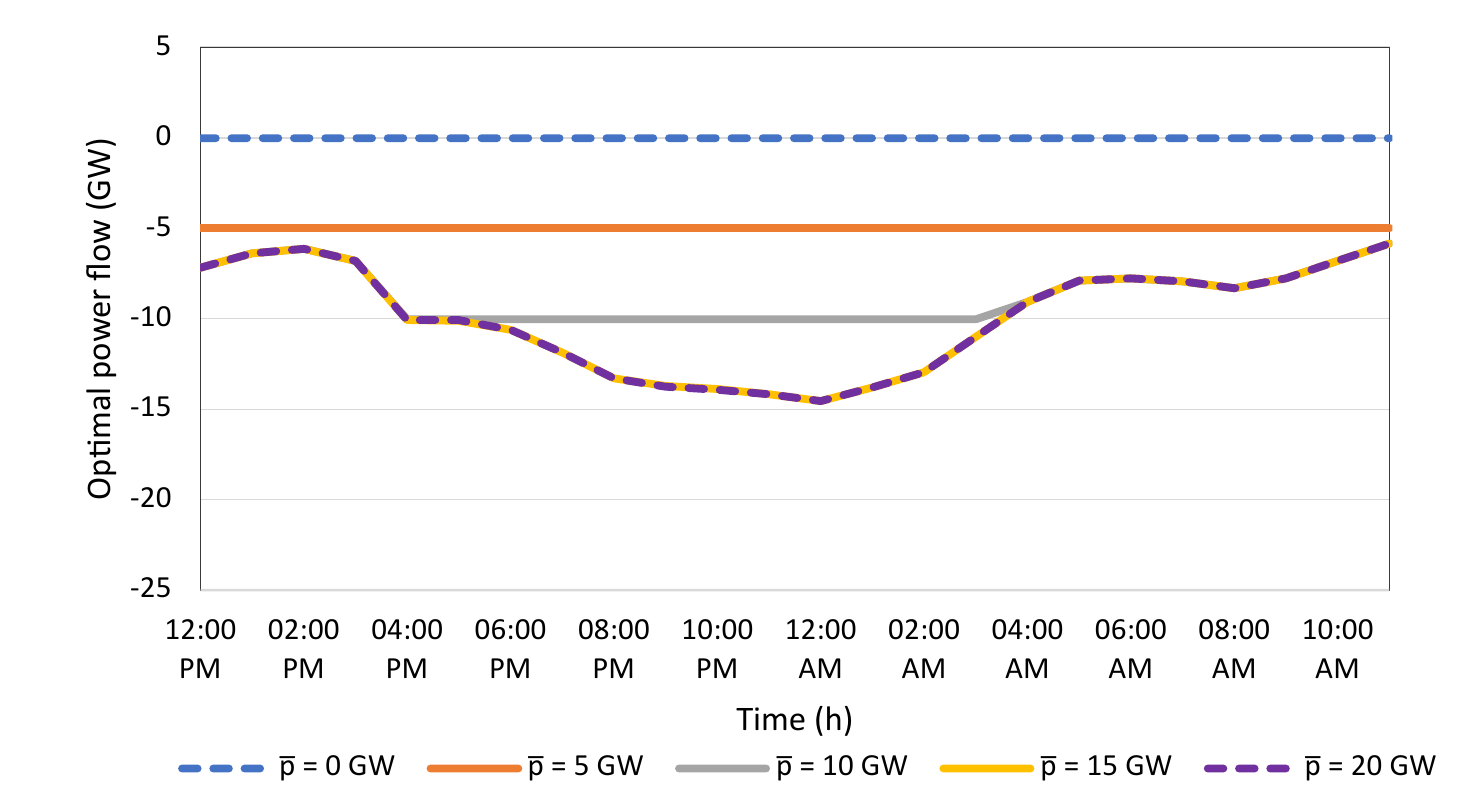}}
\vspace{-3mm}
\caption{The optimal power flow under different transmission capacities.}
\label{PowerflowTwoArea}
\end{figure}

Next, based on these methods, we discuss the impact of transmission congestion on the marginal cost at each area and on the total generation cost. Fig.~\ref{MarginalPrices} compares the electricity price profiles of the two areas. It can be seen that the marginal prices are different at each area and generally flattened by the flexible demand coordination. The price signals of the two areas converge to the same profile with the increase of power capacity so that the transmission line is no longer congested, i.e. for $\overline{p} = 15$ GW and $\overline{p} = 20$ GW. In particular, since the marginal generation costs at area 2 are higher than those of area 1, the power flow is delivered from area 1 to area 2 which is also verified by Fig.~\ref{PowerflowTwoArea}.

The optimal power flow $p(t)$ between the two areas is shown in Fig.~\ref{PowerflowTwoArea}. Negative values indicate that the power is flowing from area 1 to area 2. When the transmission capacity is no more than $5$ GW, the line is congested over 24 h for energy arbitrage. It is noted that the additional required transmission capacity to avoid congestion is around $15$ GW.

\begin{table}[H]
\centering
\vspace{-5mm}
\caption{Total generation costs (million pounds) under different Transmission Line capacity}
\scalebox{1}{
\begin{tabular}{@{}ccccccc@{}}
\toprule
\multicolumn{1}{c}{$\overline{p}$ (GW)}  & 0 & 5 & 10 & 15 & 20 \\ 
\midrule
\multicolumn{1}{c|}{MM1} & 771.69  & 718.37 & 698.14 & 696.18 & 696.18 \\
\multicolumn{1}{c|}{MM3} & 771.69 & 718.37 & 698.12 & 696.05 & 696.05  \\
\multicolumn{1}{c|}{MM4} & 771.62 & 718.29 & 698.10 & 696.18 & 696.18  \\
\bottomrule
\end{tabular}}
\label{table:TransmissionLineParameters}
\end{table}

The optimal costs using approaches MM1 - MM4 are collected in Table.~\ref{table:TransmissionLineParameters}. Since the two methods MM1 and MM2 are equivalent, only M1 is displayed in the remaining discussion. From Table.~\ref{table:TransmissionLineParameters}, we can observe that the increase of transmission capacity reduces the operational cost of the two-area system. Moreover, it is interesting to note that MM3 (MM4, respectively) achieve the same optimal cost of
MM1 for values of $\overline{p}$ such that transmission constraints are always active  (respectively never active) over the considered time window. When the transmission line is partially congested during the day, both MM3 and MM4 have small gaps with respect to the optimal cost.

\begin{figure*}[htbp!]
\begin{multicols}{2}
    \includegraphics[width=\linewidth]{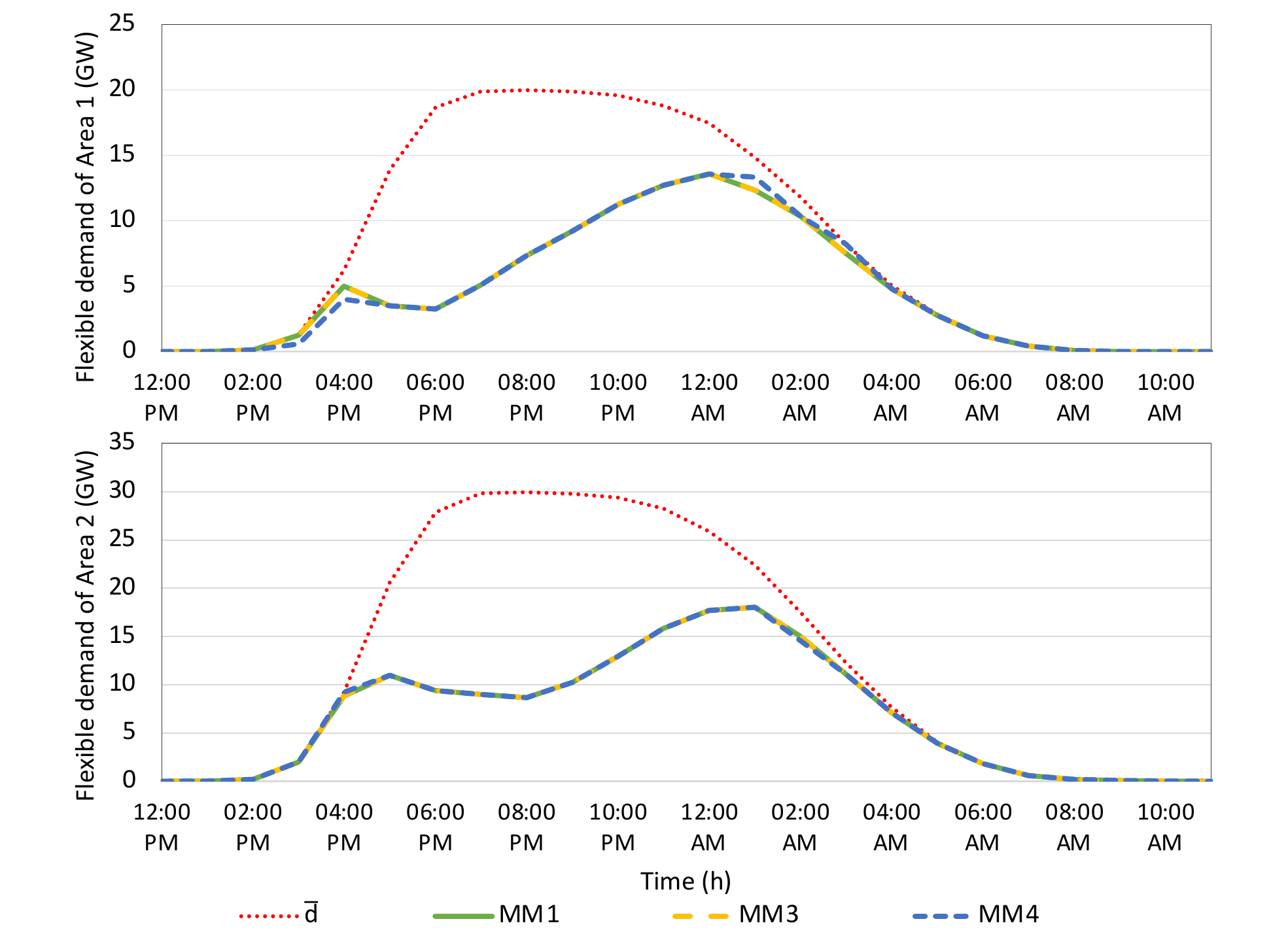}\par 
    \includegraphics[width=\linewidth]{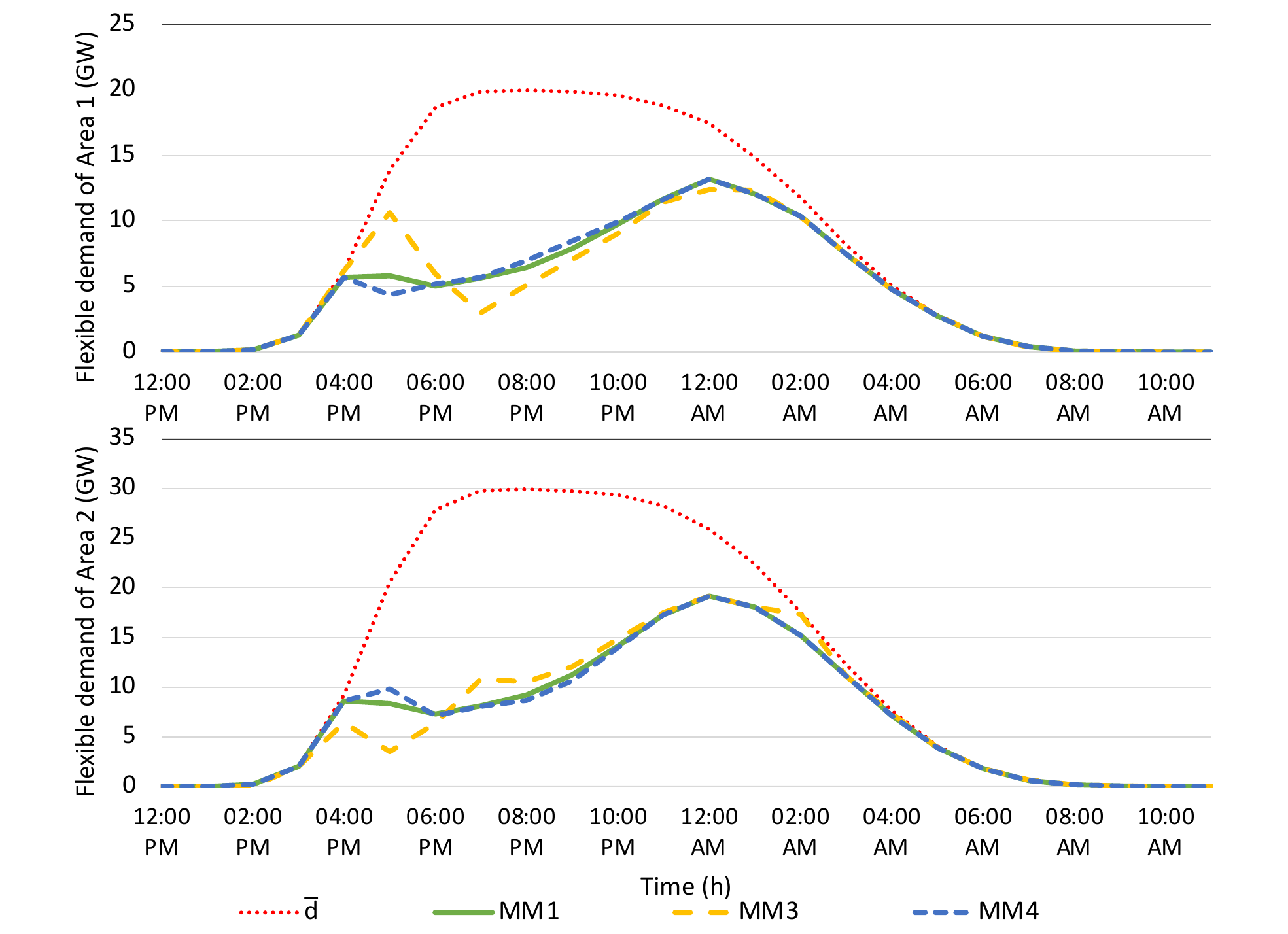}\par 
    \end{multicols}
\begin{multicols}{2}
    \includegraphics[width=\linewidth]{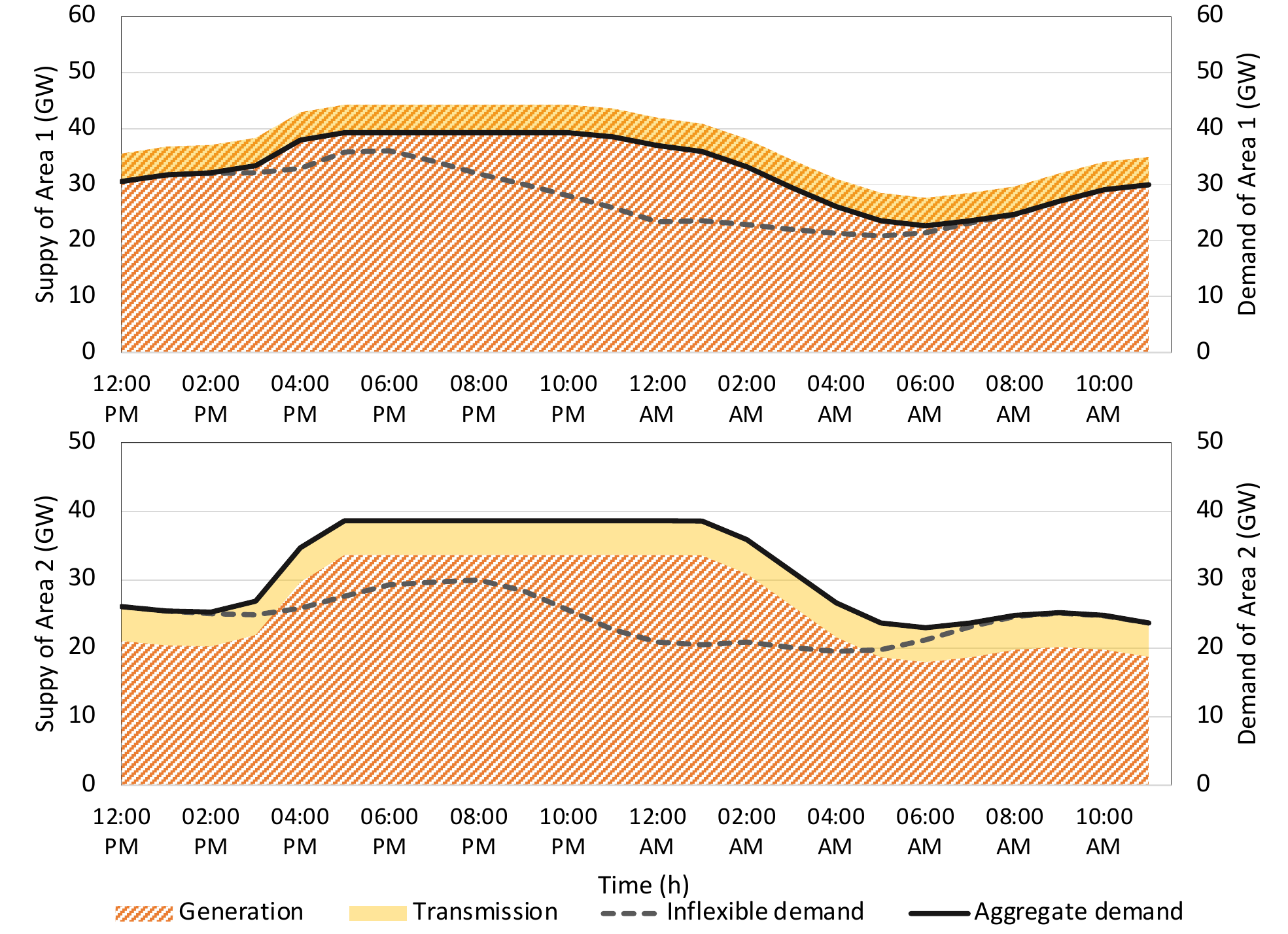}\par
    \includegraphics[width=\linewidth]{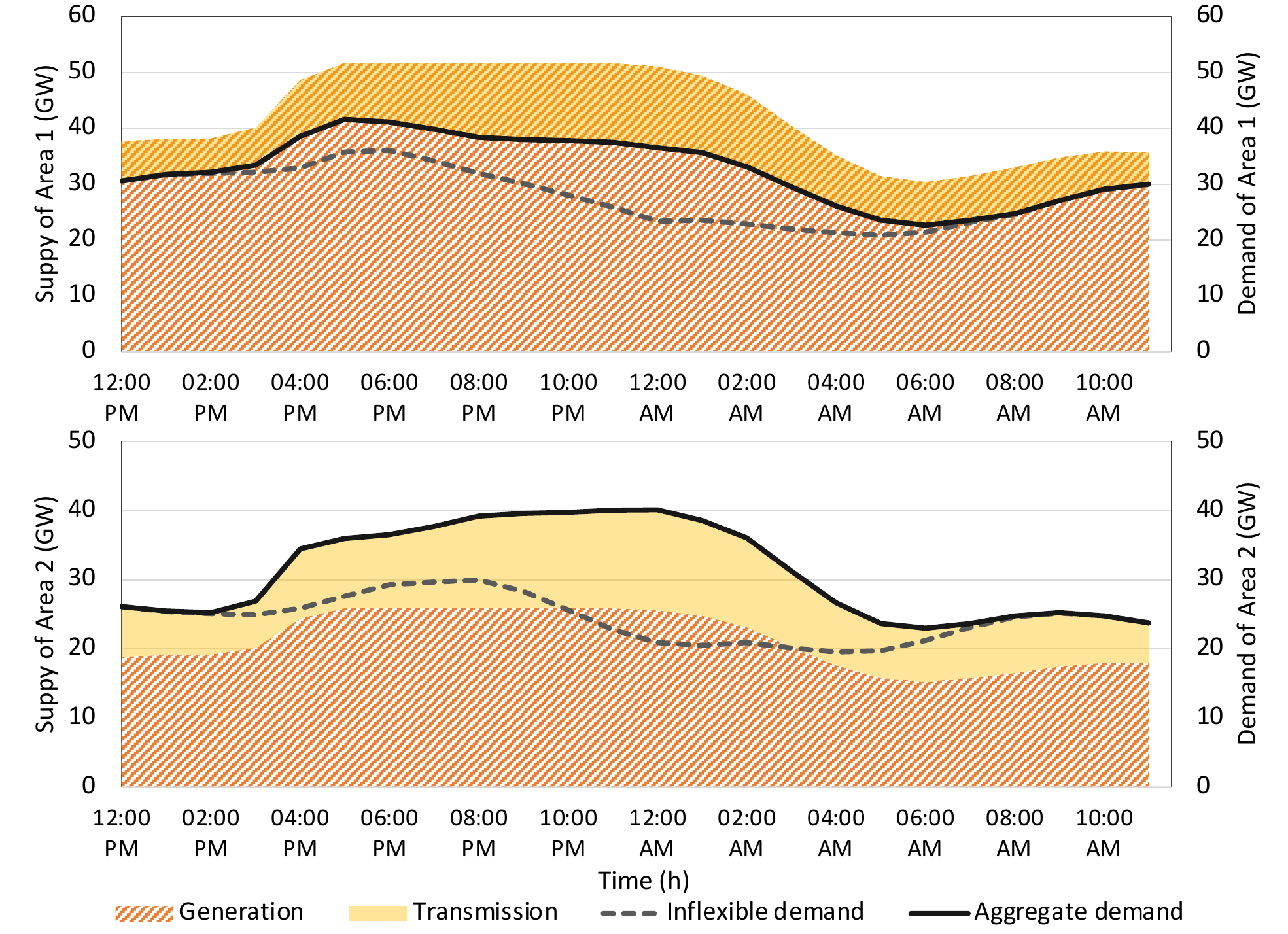}\par
\end{multicols}
\vspace{-3mm}
\caption{Top row: Aggregate demand profiles when the transmission line capacity is $5$ GW (left) and $15$ GW (right). Bottom row: Power supply and demand sources when the transmission line capacity is $5$ GW (left) and $15$ GW (right).}
\label{FlexibleDemandAndPowerBalance}
\end{figure*}

The flexible demand profiles of the two areas  are compared in Fig.~\ref{FlexibleDemandAndPowerBalance}. Note that the profiles computed by using the selection method MM3 are the same as the optimal solutions obtained by MM1, when the transmission line is congested over 24 h. However, this method may fail to obtain the optimal demand profiles if the power transferred across the transmission line is not  constant in time. In these cases, the cost is lower than the minimum cost of MM1 which indicates that the aggregate demand profiles using MM3 are unfeasible. Concerning the method MM4, which is essentially regarding the two fleets as a single merged fleet, it achieves the optimal generation cost provided that there is no congestion on the transmission line. This ensures that the total demand of the merged single fleet is feasible and optimal, but it does not necessarily imply that the demand profile in each area are feasible. To obtain the optimal aggregate demand for each fleet in the two areas, a dispatching policy of the optimal total demand to individual storage unit in the merged population is needed. Moreover, the power profiles in two areas are also included in Fig.~\ref{FlexibleDemandAndPowerBalance}, showing that the supply and demand are balanced. Compared to the results of $\overline{p} = 5$ GW, high transmission capacity allows more generation from area 1 and diverse flexible aggregate demand.

\section{Conclusion}\label{sec:conclusion}
The paper develops new methods for optimal management of heterogeneous storage fleets subject to availability constraints, by representing their aggregate behaviour through a small number of linear inequality constraints. This can be successfully applied to Unit Commitment problems where the goal is to schedule the fleet's power absorption (or, respectively power release) in order to minimize generation costs.
 The novel formulation involves a number of constraints for the fleet that scales linearly in the length $T$ of the considered time horizon and is independent of the size $N$ of the fleet. This is remarkable, as it shows that, under suitable technical conditions, a heterogeneous fleet of any size can effectively be treated as a single storage unit. The approach is compared with formulations involving a detailed description of the fleet (involving roughly $2 N T$ constraints) or aggregated descriptions involving approximately $2^T$ constraints. Numerical experiments confirm the theoretical results, showing that optimal performance can be achieved with a drastic reduction in the cost of the optimisation time. Since constraints selection is affected by a combinatorial complexity growth in the size $T$ of the considered window, we additionally propose two greedy selection procedures which are shown, through extensive simulations, to perform the correct constraint selection in most considered cases. When the choice of constraints is not optimal, this leads to optimal power schedules which are unfeasible for the original fleet (due to the reduced number of constraints considered), however numerical experiment have shown in all such cases a small optimality gap.
 The method is then tested on a two area system, where each area has its own fleet and ability to schedule it in a coordinated fashion. Transmission capacity constraints are also assumed among the two areas, so that marginal price differentials may be experienced at times when capacity constraints are active.
 While the theoretical analysis cannot be applied to this case, we show how the method yields promising results. The constraints reduction is performed by either considering the two fleets as a single fleet (i.e. ignoring capacity constraints) or, ignoring the existence of a transmission line, by carrying out the reduction of each fleet individually.

\appendix 
\section{Technical Proofs}
\subsection{Proof of Lemma \ref{optimalagent}}
Let $\mathcal{D}_a = \{ D^{I}(t)+ d^*(t):   t \in \mathcal{T} \}$. Of course,
$\mathcal{D}_a$ is a finite set. For any $j \in \mathcal{N}$ we claim that 
$d_j \in \mathcal{D}_a$ exists such that (\ref{lemmaclaim}) holds.
Let $\underline{d}_j$ be defined as:
\[  \underline{d}_j := \sup_{d \in \mathcal{D}_a}   \{ d : u_j^*(t)= \bar{P}_j \; \forall \, t: D^{I}(t) + d^*(t) \leq d    \},  \]
where $\underline{d}_j = - \infty$ if the supremum is taken over an empty set.
Similarly, we define $\bar{d}_j$ as:
\[  \bar{d}_j := \min_{d \in \mathcal{D}_a}   \{ d : u_j^*(t)= 0 \; \forall \, t: D^{I}(t) + d^*(t) \geq d    \},  \]
where $\bar{d}_j = + \infty$ if the infimum is taken over an empty set.
By construction $\underline{d}_j < \bar{d}_j$. 
We need to show that there exists at most a single element $d \in \mathcal{D}_a$ with the property that $\underline{d}_j < d < \bar{d}_j$. If the claim is true, indeed, (\ref{lemmaclaim}) holds by choosing $d_j :=d$.
By contradiction, should this not be the case, there would exist
$t_1$ and $t_2$ in $\mathcal{T}$, with $D^{I}(t_1)+ d^*(t_1) < D^{I}(t_2) + d^*(t_2)$ such that	$u_j^*(t_1) < \bar{P}_j$ and $u_j^*(t_2)>0$.  
This, however, contradicts optimality of $u_j^*(t)$ as for any sufficiently small $\delta>0$ the policy obtained by taking
$\tilde{u}_j (t_1) = u_j^*(t_1) + \delta$, $\tilde{u}_j(t_2)
= u_j^*(t_2) - \delta$ and $\tilde{u}_j (t) = u^*_j(t)$ for all other $t \in \mathcal{T}$, is feasible and of strictly lower cost. 
\subsection{Proof of Lemma \ref{activeconstraints} }
Let $u_j^*(t)$, for $j \in \mathcal{N}$ be any optimal solution of (\ref{original}). As already remarked this fulfills
$\sum_{j \in \mathcal{N}} u_j^*(t) = d^*(t)$.
Hence, by virtue of Lemma \ref{optimalagent}
\begin{equation*}
\begin{aligned}
\sum_{t \in W_d} d^*(t) &= \sum_{t \in W_d} \sum_{j \in \mathcal{N}}
u_j^*(t)\\
&= \sum_{j \in \mathcal{N}} \sum_{t \in W_d} u_j^*(t) \\
&= \sum_{j: d \leq d_j} \sum_{t \in W_d} u_j^*(t)   +
\sum_{j: d > d_j} \sum_{t \in W_d} u_j^*(t)\\
&=   \sum_{j: d \leq d_j} \textrm{card} ( \mathcal{A}_j \cap W_d ) \bar{P}_j + \sum_{j: d > d_j} \bar{E}_j \\ & \geq \sum_{j \in \mathcal{N}}
	\min \{ \textrm{card} ( \mathcal{A}_j \cap W_d ), \bar{E}_j/ \bar{P}_j \} \bar{P}_j.
\end{aligned}
\end{equation*}
% \[ \sum_{t \in W_d} d^*(t) = \sum_{t \in W_d} \sum_{j \in \mathcal{N}}
% u_j^*(t)  \]
% \[ \qquad = \sum_{j \in \mathcal{N}} \sum_{t \in W_d} u_j^*(t)    \]
% \[ \qquad  = \sum_{j: d \leq d_j} \sum_{t \in W_d} u_j^*(t)   +
% \sum_{j: d> d_j} \sum_{t \in W_d} u_j^*(t)   \]	
% \[ \qquad =   \sum_{j: d \leq d_j} \textrm{card} ( \mathcal{A}_j \cap W_d ) \bar{P}_j + \sum_{j: d> d_j} \bar{E}_j \geq \sum_{j \in \mathcal{N}}
% 	\min \{ \textrm{card} ( \mathcal{A}_j \cap W_d ), \bar{E}_j/ \bar{P}_j \} \bar{P}_j.  \]
Hence, since $d^*(t)$ is feasible, $W_d$ is the support of an active constraint.	\\
\subsection{Proof of Lemma \ref{achieved}}
Let $d^*(t)$ be the optimal aggregated demand for problem (\ref{modified}). Then, by Lemma \ref{activeconstraints} the following holds
for any $W \subset \mathcal{T}$ of the same cardinality as $W_d$:
\[ \frac{\sum_{t \in W_d } D^{I}(t) +
	\sum_{j \in \mathcal{N}} \min \{ \textrm{card} ( \mathcal{A}_j \cap W_d ), \bar{E}_j / \bar{P}_j \} \bar{P}_j}{\textrm{card} ( W_d)}   \]
\[    =   \frac{\sum_{t \in W_d} D^{I}(t) + d^*(t) }{\textrm{card}(W_d)}
\leq \frac{\sum_{t \in W} D^{I}(t) + d^*(t) }{ \textrm{card} (W) }
  \leq \bar{D} ( W).\]
Hence, 
\[ \bar{D} (W_{d}) \leq  \min_{W \subset \mathcal{T}: \textrm{card} (W)= \textrm{card}(W_d) }   \bar{D} (W).\]
The final claim follows by remarking that the inequality
\[ \frac{\sum_{t \in W_d} D^{I}(t) + d^*(t) }{\textrm{card}(W_d)}
\leq \frac{\sum_{t \in W} D^{I}(t) + d^*(t) }{ \textrm{card} (W) } \]
is strict whenever $W \neq W_d$, (this is because $W$ must include at least a time instant with strictly greater $D^{I}(t) + d^*(t)$ than the corresponding values for $W_d$, by definition of $W_d$). 
\subsection{Proof of Lemma \ref{wdenough}}
Let $\tilde{d}^*(t)$ be the optimal solution of (\ref{modified2}). 
We claim that for any $k \in \{1, \ldots, Q-1 \}$ the signal $D^{I}(t) + \tilde{d}^*(t)$ is constant on $W_{k+1} \backslash W_{k}$.
This claim is trivial if $W_{k+1} \backslash W_k$ has cardinality $1$.
If, cardinality is bigger than one we see arguing by contradiction that a sufficiently small power swap from the time instant where the maximum is achieved towards the time instant where the minimum is achieved would still fulfil all constraints and imply a net cost reduction, thus violating optimality of $\tilde{d}^*(t)$.
Moreover, $d^*(t) = \tilde{d}^*(t)$ for all $t$. We prove the result $W_k$ by induction showing that,
\[ \sum_{t \in W_k} \tilde{d}^*(t) = \sum_{t \in W_k} d^*(t), \]
viz. the constraint of support $W_k$ is active.
For $k=1$ notice that
$\sum_{t \in W_1} \tilde{d}^*(t) < \sum_{t \in W_1} d^*(t)$
implies that $D^{I}(t) + \tilde{d}^*(t) > D^{I}(t) + d^*(t)$ for some $t \notin W_1$. The latter value, is in turn bigger than $D^{I}(t) +d(t)$ for $t \in W_1$. Hence a power swap from such $t$ into all time instants in $W_1$ (equally distributed) would still preserve constraints and yield an overall cost reduction, violating optimality of $\tilde{d}^*(t)$. 
In particular, then $\tilde{d}^*(t)= d^*(t)$ for $t \in W_1$.
By induction, a similar argument shows that
 $\sum_{t \in W_k} \tilde{d}^*(t) = \sum_{t \in W_k} d^*(t)$. Assume by hypothesis that the result is true for $k-1$. If, by contradiction $\sum_{t \in W_k} \tilde{d}^*(t) < \sum_{t \in W_k} d^*(t)$,
 then $\sum_{t \in W_{k} \backslash W_{k-1} } \tilde{d}^*(t)
 < \sum_{t \in W_k \backslash W_{k-1} } d^*(t)$ implies that 
 $D^{I}(t) + \tilde{d}^*(t) > D^{I}(t) + d^*(t)$ for some $t \notin W_k$.
The latter value, is in turn bigger than $D^{I}(t) +d(t)$ for $t \in W_k$. Hence a power swap from such $t$ into all time instants in $W_k \backslash W_{k-1}$ (equally distributed) would still preserve constraints and yield an overall cost reduction, violating optimality of $\tilde{d}^*(t)$.
As a consequence $\tilde{d}^*(t) = d^*(t)$, for all $t \in W_k$.
This proves the claim and the Lemma. 
\subsection{Proof of Theorem \ref{wstarenough}}
Thanks to the inclusion (\ref{nested}) we see that
the optimal cost of (\ref{modified2}) is less or equal to the optimal cost of
(\ref{modified3}) (because (\ref{modified3}) is minimised over a subset of the feasible region of (\ref{modified2})).
On the other hand, again by virtue of (\ref{nested}), the optimal cost of (\ref{modified3}) is less or equal to the optimal cost of (\ref{modified}),
since (\ref{modified}) is minimised over a subset of the feasible region of (\ref{modified3}). Hence, the claim follows by Lemma \ref{wdenough}.

\begin{IEEEbiography}[{\includegraphics[width=1in,height=1.25in,clip,keepaspectratio]{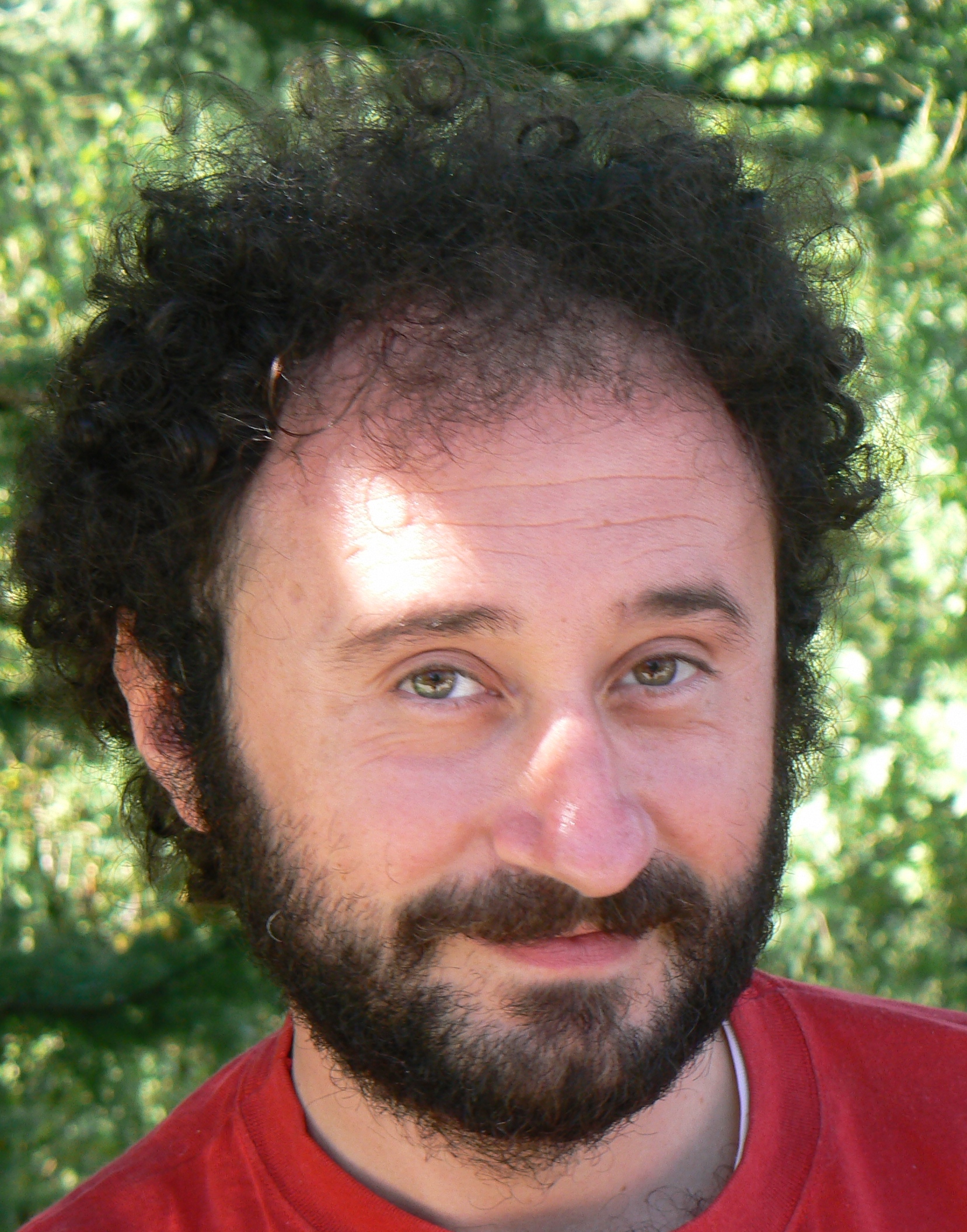}}]{David Angeli} was born in Siena, Italy, in 1971. He received the B.S. degree in Computer Science Engineering and the Ph.D. in Control Theory from University of Florence, Florence, Italy, in 1996 and 2000, respectively.

 Since 2000 he was an Assistant and Associate Professor (2005) with the Department of Systems and Computer Science, University of Florence. In 2007 he was a visiting Professor with I.N.R.I.A de Rocquencourt, Paris, France, and since 2008, he joined as a Senior Lecturer the Department of Electrical and Electronic Engineering of Imperial College London, where he is currently a Reader in Nonlinear Control and the Director of the MSc in Control Systems. He has been a Fellow of the IEEE since 2015.

 He has served as an Associate Editor for IEEE Transactions in Automatic Control and Automatica. Overall he authored more than 100 journal papers in the areas of Stability of nonlinear systems, Control of constrained systems (MPC), Chemical Reaction Networks theory and Smart Grids.
\end{IEEEbiography}

\vspace{-3mm}

\begin{IEEEbiography}[{\includegraphics[width=1in,height=1.25in,clip,keepaspectratio]{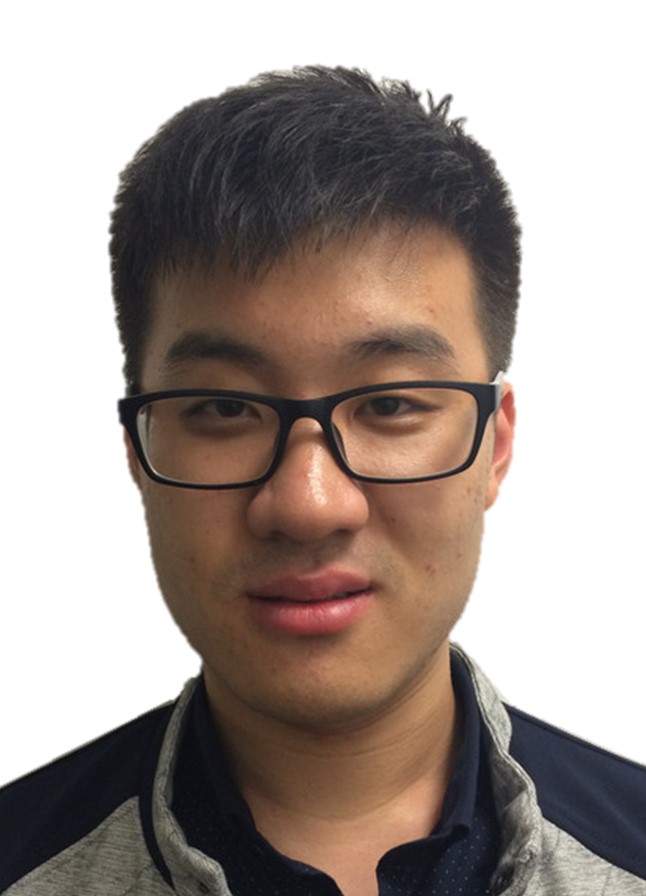}}]{Zihang Dong} received B.Eng degree in Electrical and Electronics Engineering in 2014 from the University of Liverpool, UK. He received M.Sc. Degree and the Ph.D. degree in Control Theory from Imperial College London, UK, in 2015 and 2019, respectively. He is currently a research associate at the Department of Electrical and Electronic Engineering of Imperial College London. His research interests include robust economic model predictive control, flexible demand response and the design of control strategies for the smart grids.
\end{IEEEbiography}

\vspace{-3mm}

\begin{IEEEbiography}[{\includegraphics[width=1in,height=1.25in,clip,keepaspectratio]{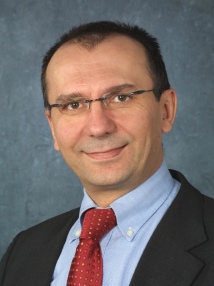}}]{Goran Strbac} is a Professor of Electrical Energy Systems at Imperial College London, U.K. His current research is focused on optimisation of energy systems operation and investment, energy infrastructure reliability and future energy markets.
\end{IEEEbiography}

\end{document}